\documentclass[10pt,reqno,fullpage]{amsart}
\usepackage{amsmath,amsthm,amssymb}%,psfig}
\usepackage{pifont
}
\theoremstyle{plain}
\newtheorem{theorem}{Theorem}[section]

\newtheorem{proposition}{Proposition}[section]

\newtheorem{definition}{Definition}[section]

\theoremstyle{remark}
\newtheorem{remark}{Remark}[section]

\newfont{\menge}{dsrom10}

\newcommand{\R}{\text{\menge R}}
\newcommand{\eps}{\varepsilon}
\numberwithin{equation}{section}

\begin{document}

\title[Hyperbolic Conservation Laws and Hydrodynamic Limit for Particle Systems ]
{Hyperbolic Conservation Laws with Discontinuous Fluxes and
Hydrodynamic Limit
  for Particle Systems }
\author{Gui-Qiang Chen}
\address{G.-Q. Chen, Department of Mathematics, Northwestern University, 2033 Sheridan Road,
         Evanston, IL 60208-2730, USA}
\email{gqchen@math.northwestern.edu}
%\thanks{}
\author{Nadine Even}
\address{N. Even, Department of
Mathematics, University of W\"urzburg, Am Hubland, D--97074
W\"u{}rzburg, Germany} \email{even@mathematik.uni-wuerzburg.de}
%\email{nadineeven@gmx.de}
\author{Christian Klingenberg}
\address{C. Klingenberg, Department of
Mathematics, University of W\"urzburg, Am Hubland, D--97074
W\"u{}rzburg, Germany}
\email{klingenberg@mathematik.uni-wuerzburg.de} \subjclass{ 35L65,
60K35, 82C22} \keywords{hyperbolic conservation laws, discontinuous
flux functions, measure-valued, entropy solutions, entropy
conditions, uniqueness, hydrodynamic limits, microscopic, particle
systems, zero range process, process of misanthropes, compactness
framework}

\begin{abstract}
We study the following class of scalar hyperbolic conservation laws
with discontinuous fluxes:
\begin{equation}\label{0.1}
\partial_t\rho+\partial_xF(x,\rho)=0.
\end{equation}
The main feature of such a conservation law is the discontinuity of
the flux function in the space variable $x$. Kruzkov's approach
%in\cite{kruzkov}
for the $L^1$-contraction does not apply since it requires the
Lipschitz continuity of the flux function; and entropy solutions
even for the Riemann problem are not unique under the classical
entropy conditions. On the other hand, it is known that, in
statistical mechanics, some microscopic interacting particle systems
with discontinuous speed parameter $\lambda(x)$, in the hydrodynamic
limit, formally lead to scalar hyperbolic conservation laws with
discontinuous fluxes of the form:
 \begin{equation}\label{0.2}
\partial_t\rho+\partial_x\left(\lambda(x)h(\rho)\right)=0.
\end{equation}
The natural question arises which entropy solutions the hydrodynamic
limit selects, thereby leading to a suitable, physical relevant
notion of entropy solutions of this class of conservation laws. This
paper is a first step and provides an answer to this question for a
family of discontinuous flux functions. In particular, we identify
the entropy condition for \eqref{0.1} and proceed to show the
well-posedness by combining our existence result with a uniqueness
result of Audusse-Perthame (2005) for the family of flux functions;
we establish a compactness framework for the hydrodynamic limit of
large particle systems and the convergence of other approximate
solutions to \eqref{0.1}, which is based on the notion and reduction
of measure-valued entropy solutions; and we finally establish the
hydrodynamic limit for a ZRP with discontinuous speed-parameter
governed by an entropy solution to \eqref{0.2}.
\end{abstract}

\maketitle

%%%%%%%%%%%%%%%%%%%%%%%%%%%%%%
\section{Introduction}
%%%%%%%%%%%%%%%%%%%%%%%%%%%%%%

We are concerned with the following class of scalar hyperbolic
conservation laws with discontinuous fluxes:
\begin{equation}\label{1.1}
\partial_t\rho+\partial_xF(x,\rho(t,x))=0
\end{equation}
and with initial data:
\begin{equation}\label{1.1a}
\rho|_{t=0}=\rho_0(x),
\end{equation}
where $F(\cdot,\rho)$ is continuous except on a set of measure zero.

The main feature of \eqref{1.1} is the discontinuity of the flux
function in the space variable $x$. This feature causes new
important difficulties in conservation laws. Kruzkov's approach in
\cite{kruzkov} for the $L^1$-contraction does not apply; entropy
solutions even for the Riemann problem  of \eqref{1.1} are not
unique under the classical entropy conditions; several admissibility
criteria have been proposed in
\cite{ap,BaitiJenssen,GimseRisebro,KarlsenRisebroTowers,KRi} and the
references cited therein. In particular, a uniqueness theorem was
established in Baiti-Jenssen \cite{BaitiJenssen} when $F(x,\cdot)$
is monotone and Audusse-Perthame \cite{ap} for more general flux
functions that especially include non-monotone functions
$F(x,\cdot)$ in \eqref{1.1} under their notion. However, the
existence of entropy solutions for the non-monotone case under the
notion of Audusse-Perthame \cite{ap} has not been established, and
the entropy conditions proposed in the literature in general are not
equivalent.

On the other hand, in statistical mechanics, some microscopic
interacting particle systems with discontinuous speed parameter
$\lambda(x)$, in the hydrodynamic limit, formally lead to scalar
hyperbolic conservation laws with discontinuous flux of the form
 \begin{equation}\label{1.2}
\partial_t\rho+\partial_x\left(\lambda(x)h(\rho)\right)=0
\end{equation}
and with initial data \eqref{1.1a},
where $\lambda(x)$ is continuous except on a set of measure zero and
$h(\rho)$ is Lipschitz continuous. Equation \eqref{1.2} is
equivalent to the following $2\times 2$ hyperbolic system of
conservation laws:
\begin{equation}\label{system1.4}
\left\{\begin{array}{ll}
\partial_t\rho+ \partial_x(\lambda h(\rho))=0, \\
\partial_t\lambda = 0.
\end{array}\right.
\end{equation}
In particular, when $h(\rho)$ is not strictly monotone, system
\eqref{system1.4} is nonstrictly hyperbolic, one of the main
difficulties in conservation laws (cf. \cite{chen,Da}). The natural
question is which entropy solution the hydrodynamic limit selects,
thereby leading to a suitable, physical relevant notion of entropy
solutions of this class of conservation laws. This paper is a first
step and provides an answer to this question for a family of
discontinuous flux functions via an interacting particle system,
namely, the attractive zero range process (ZRP). This ZRP leads to a
conservation law of the form \eqref{1.2} with $\lambda(x)>0$ and
$h(\rho)$ being monotone in $\rho$, and its hydrodynamic limit
naturally gives rise to an entropy condition of the type described
in \cite{ap,BaitiJenssen} in the formal level.

Motivated by the hydrodynamic limit of the ZRP, in this paper, we
adopt the notion of entropy solutions for a class of conservation
laws with discontinuous flux functions, including the non-monotone
case in the sense of Audusse-Perthame \cite{ap}, and establish the
existence of such an entropy solution via the method of compensated
compactness in Section 3. This completes the well-posedness by
combining a uniqueness result established in \cite{ap} for this
class of conservation laws under the notion of entropy solutions.

In order to establish the hydrodynamic limit of large particle
systems and the convergence of other approximate solutions to
\eqref{1.1} rigorously, we establish another compactness framework
for \eqref{1.1}--\eqref{1.1a} in Section 2. This mathematical
framework is based on the notion and reduction of measure-valued
entropy solutions developed in Section 2, which is also applied for
another proof of the existence of entropy solutions for the monotone
case in Section 3.

In Section 4, we establish the hydrodynamic limit for a ZRP with
discontinuous speed-parameter $\lambda(x)$ governed by the unique
entropy solution of the Cauchy problem \eqref{1.1a}--\eqref{1.2}.

\section{Notion and Reduction of measure-valued entropy solutions}
\label{section1}

In this section, we first develop the notion of measure-valued
entropy solutions and establish their reduction to  entropy
solutions in $L^{\infty}$ (provided that they exist) of the Cauchy
problem \eqref{1.1}--\eqref{1.1a} satisfying
\begin{itemize}
\item[(H1)] $F(x,\rho)$ is continuous at all points of
  $({\text{\menge R}}\backslash\mathcal N)\times {\text{\menge R}}$ with
  $\mathcal N$ a closed set of measure zero;
\item[(H2)] $\exists$ continuous functions $f,g$ such that,
  for any $x\in{\text{\menge R}}$ and large $\rho$,
  $f(\rho)\leq|F(x,\rho)|\leq g(\rho)$ with $f(\rho)\ge 0$
  and $f(\pm\infty)=\infty$;
\item [(H3)] There exists a function $\rho_m(x)$ from ${\text{\menge R}}$ to ${\text{\menge R}}$
  and a constant $M_0$ such that,
 for $x\in{\text{\menge R}}\backslash\mathcal N$, $F(x,\rho)$ is a locally Lipschitz,
  one to one function from $(-\infty,\rho_m]$ and $[\rho_m, \infty)$ to
  $[M_0,\infty)$ (or $(-\infty, M_0]$) with $F(x,\rho_m(x))=M_0$;
\end{itemize}
or
\begin{itemize}
\item[(H3')] For $x\in {\text{\menge R}}\backslash\mathcal N$,
$F(x,\cdot)$ is a locally Lipschitz, one to one function from
${\text{\menge R}}$ to ${\text{\menge R}}$.
\end{itemize}

\smallskip One example of the flux function satisfying
(H1)--(H2) and (H3) or (H3') is
\begin{equation}
F(x,\rho)=\lambda(x)h(\rho),
\end{equation}
where $\lambda(x)$ is continuous in $x\in\text{\menge R}$ with
$0<\lambda_1\le \lambda(x)\le \lambda_2<\infty$ for some constants
$\lambda_1$ and $\lambda_2$, except on a closed set $\mathcal N$ of
measure zero, $h(\rho)$ is locally Lipschitz and is either monotone
or convex (or concave) with $h(\rho_m)=0$ for some $\rho_m$ in which
case $M_0=0$.

It is easy to check that, if the flux function $F(x,\rho)$ satisfies
(H1)--(H3), then, for any constant $\alpha\in [M_0,\infty)$ (or
$\alpha\in (-\infty, M_0]$), there are two steady-state solutions
$m_{\alpha}^+$ from ${\text{\menge R}}$ to $[\rho_m(x),\infty)$ and
$m_{\alpha}^-$ from ${\text{\menge R}}$ to $(-\infty,\rho_m(x)]$ of
\eqref{1.1} such that
\begin{equation}\label{m1}
F(x,m^{\pm}_{\alpha}(x))=\alpha.
\end{equation}
In the case (H1)--(H2) and (H3'), $m_\alpha^+(x)=m_\alpha^-(x)$
which is even simpler.

%%%%%%%%%%%%%%%%%%%%%%%%%%%%%%%%%%%%%%%
\subsection{Notion of measure-valued entropy solutions}
%%%%%%%%%%%%%%%%%%%%%%%%%%%%%%%%%%%%%%%%%%%%

First, the notion of entropy solutions in $L^\infty$ introduced in
Audusse-Perthame \cite{ap} and Baiti-Jenssen \cite{BaitiJenssen} can
be further formulated into the following.

\begin{definition}[Notion of entropy solutions in
$L^\infty$]\label{weakap} We say that an $L^{\infty}$ function
$\rho: \R_+^2:={\text{\menge R}}_+\times{\text{\menge
R}}\mapsto{\text{\menge R}}$ is an entropy solution of
\eqref{1.1}--\eqref{1.1a} provided that, for each $\alpha\in
[M_0,\infty)$ (or $\alpha\in (-\infty, M_0]$) and the corresponding
two steady-state solutions $m_{\alpha}^{\pm}(x)$ of \eqref{1.1},
\begin{eqnarray}
&&\int\Big({|\rho(t,x)-m^{\pm}_{\alpha}(x)|}\,\partial_tJ\,
+{\text{sgn}}(\rho(t,x)-m^{\pm}_{\alpha}(x))\big(F(x,\rho(t,x))-\alpha\big)\,\partial_x
J\Big)\, dtdx
\nonumber\\
&&\qquad +\int \, |\rho_0(x)-m^{\pm}_{\alpha}(x)| J(0,x)\, dx \geq0
\label{es-2}
\end{eqnarray}
for any test function $J: \R^2_+\mapsto \R_+$.
\end{definition}

\smallskip It is easy to see that any entropy solution is a
weak solution of \eqref{1.1}--\eqref{1.1a} by choosing $\alpha$ such
that $m_\alpha^+(x)\ge \|\rho\|_{L^\infty}$ and $m_\alpha^-(x)\le
-\|\rho\|_{L^\infty}$, respectively, for a.e. $x\in \R$.

{}From the uniqueness argument in Audusse-Perthame \cite{ap} (also
see \cite{ChenRascle}), one can deduce that, for any $L>0$,
\begin{equation}\label{initial-data}
\lim_{t\to 0}\int_{|x|\le L}|\rho(t,x)-\rho_0(x)|\, dx =0.
\end{equation}

\smallskip Following the notion of entropy solutions, we
introduce the corresponding notion of measure-valued entropy
solutions. We denote by $\mathcal P({\text{\menge R}})$ the set of
probability measures on ${\text{\menge R}}$.

\begin{definition}[Notion of measure-valued entropy
solutions]\label{mvap}
We say that a measurable map
$$
\pi:\, \R_+^2\rightarrow{\mathcal P({\text{\menge R}})}
$$
is a measure-valued entropy solution of \eqref{1.1}--\eqref{1.1a}
provided that $\langle \pi_{0,x}; k\rangle=\rho_0(x)$ for a.e. $x\in
\R$ and, for each $\alpha\in[M_0, \infty)$ (or $\alpha\in (-\infty,
M_0]$) and the corresponding two steady-state solutions
$m_\alpha^\pm(x)$ of \eqref{1.1},
\begin{eqnarray}
&&\int\left(\langle\pi_{t,x};
 |k-m^{\pm}_{\alpha}(x)|\rangle \partial_tJ
+\langle\pi_{t,x}; {\text{sgn}}(k-m^{\pm}_{\alpha}(x))
\left(F(x,k)-\alpha\right)\rangle \partial_x J\right)dxdt\nonumber\\
&&\qquad +\int |\rho_0(x)-m^{\pm}_{\alpha}(x)|\, J(0,x)\, dx \ge 0
\label{mves-2}
\end{eqnarray}
for any test function $J: \R_+^2\mapsto{\text{\menge R}_+}$.
\end{definition}

If a measure-valued entropy solution $\pi_{t,x}(k)$ is a Dirac mass
with the associated profile $\rho(t,x)$, i.e.
$\pi_{t,x}(k)=\delta_{\rho(t,x)}(k)$, then $\rho(t,x)$ is an entropy
solution of \eqref{1.1}--\eqref{1.1a}, which is unique as shown in
\cite{ap}.

Note that, when the flux function $F(x,\rho)$ is locally Lipschitz
in both variables $(x,\rho)$,  one can use the Kruzkov entropy
inequality, instead of (\ref{mves-2}),  to formulate the following
notion of measure-valued solutions:
\begin{equation}\label{Kr}
\partial_t\left\langle\pi_{t,x}; \left|k-c\right|\right\rangle
+\partial_x\left\langle\pi_{t,x};
{\text{sgn}}(k-c)\left(F(x,k)-F(x,c)\right)\right\rangle
+\left\langle\pi_{t,x};
{\text{sgn}}(k-c)\partial_xF(x,c)\right\rangle \leq0
\end{equation}
in the sense of distributions and to establish their reduction as in
DiPerna \cite{diperna}. One of the new features in our formulation
\eqref{mves-2} in Definition 2.2 is that the constant $c$ in
\eqref{Kr} is replaced by the steady-state solutions
$m^\pm_\alpha(x)$ so that the additional third term in \eqref{Kr}
vanishes, as in \cite{ap,BaitiJenssen}, and thereby allows the
discontinuity of the flux functions on a closed set of measure zero
for measure-valued entropy solutions.

\subsection{Reduction of measure-valued entropy solutions}

In this section we first establish the reduction of measure-valued
entropy solutions of \eqref{1.1}--\eqref{1.1a} and prove that any
measure-valued entropy solution $\pi_{t,x}(k)$ in the sense of
Definition 2.2 is the Dirac solution such that the associated
profile $\rho(t,x)$ is an entropy solution in the sense of
Definition 2.1. That is, our goal is to establish that, when
$\pi_{0,x}(k)=\delta_{\rho_0(x)}(k)$,
\begin{equation}\label{reduction}
\pi_{t,x}(k)=\delta_{\rho(t,x)}(k),
\end{equation}
where $\rho:\,\R_+^2\rightarrow{\text{\menge R}}$ is the unique
entropy solution determined by \eqref{es-2}. The reduction proof is
achieved by two theorems. We start with the first theorem.

{\theorem \label{dip} Assume $\rho:\, \R_+^2\rightarrow{\text{\menge
R}}$ is the unique entropy solution of \eqref{1.1}--\eqref{1.1a}
with initial data $\rho_0\in L^{\infty}({\text{\menge R}})$. Assume
that there exists a measure-valued entropy solution $\pi:\,
\R_+^2\rightarrow\mathcal P({\text{\menge R}})$ of \eqref{1.1} in
the sense of Definition {\rm \ref{mvap}} with $\pi_{t,x}$ having a
fixed compact support for a.e. $(t,x)$ and
$\pi_{0,x}(k)=\delta_{\rho_0(x)}(k)$ for a.e. $x\in \R$. Then
\begin{equation}
\int \big(\left\langle\pi_{t,x};
 \left|k-\rho(t,x)\right|\right\rangle\partial_tJ
 +\left\langle\pi_{t,x}; {\text
{sgn}}(k-\rho(t,x))(F(x,k)-F(x, \rho(t,x)))\right\rangle \partial_x
J\big)\, dxdt\ge 0 \label{E}
\end{equation}
for any test function $J: \,\R_+^2\mapsto{\text{\menge R}_+}$. }

\begin{proof} The proof is divided into six steps.

{\it Step 1.} We first rewrite
$$
E:=\partial_t\big\langle\pi_{t,x};|k-\rho(t,x)|\big\rangle
+\partial_x\big\langle\pi_{t,x};{\text{sgn}}\left(k-\rho(t,x)\right)
\big(F(x,k)-F(x,\rho(t,x))\big)\big\rangle
$$
in the entropy inequality \eqref{E}. We notice the following:
\begin{itemize}
\item Under the assumption (H3'), $F(x,\rho(t,x))$ is continuous in $x$ a.e.
Then we can define a function $\tilde{\rho}(s,y,x)$ for a.e.
$(s,y,x)\in \R_+\times\R^2$ such that, for fixed $(s,y)$,
\begin{eqnarray}\label{52}
F(x,\tilde{\rho}(s,y,x)):=F(x,m_{F(y,\rho(s,y))}(x))=F(y,\rho(s,y)),
\end{eqnarray}
where the last equality follows from \eqref{m1}. Thus, we define
$$
\beta(s,y):=F(y,\rho(s,y)) \quad\mbox{so that}\quad
\tilde{\rho}(s,y,x)=m_{\beta(s,y)}(x).
$$
\item For the case (H3), we define the sign of the difference between
the tilda function and $\rho_m(y)$ to be the same as the sign of the
corresponding solution. Since $\rho_m(y)$ is the minimum (or
maximum) point of the flux function with $F(y,\rho_m(y))=M_0$, then,
for
\begin{eqnarray}\label{tilde-rho}
\tilde{\rho}(s,y,x):=m^{+}_{\beta(s,y)}(x){\text{sgn}}_+(\rho(s,y)-\rho_m(y))
+m^-_{\beta(s,y)}(x){\text{sgn}}_-(\rho(s,y)-\rho_m(y)),
\end{eqnarray}
we have as in \eqref{52}
$$
\begin{array}{ll}
&F(x,\tilde{\rho}(s,y,x)):=F(x,m^{+}_{\beta(s,y)}(x){\text{sgn}}_+(\rho(s,y)-\rho_m(y))
+m^-_{\beta(s,y)}(x){\text{sgn}}_-(\rho(s,y)-\rho_m(y))\\
&\qquad\qquad\qquad\,\,\,\,\,=F(y,\rho(s,y))=\beta(s,y).
\end{array}
$$
\end{itemize}
With these notations, we set
\begin{equation}\label{tildeE}
\tilde{E}:=\partial_t\big\langle\pi_{t,x};
\left|k-\tilde{\rho}(s,y,x)\right|\big\rangle
+\partial_x\big\langle\pi_{t,x};{\text
{sgn}}(k-\tilde{\rho}(s,y,x))\big(F(x,k)-\beta(s,y)\big)\big\rangle.
\end{equation}
Then, to obtain the
inequality $E\leq 0$, it suffices to show that $\lim_{x\rightarrow
y}\tilde E=E$.

\medskip

{\it Step 2.} We now show that
\begin{equation}\label{40}
\tilde{\rho}(s,y,x)\stackrel{x\longrightarrow
y}{\longrightarrow}\tilde{\rho}(s,y,y)=\rho(s,y)\qquad \mbox{for
a.e.}\,\, (s,y)\in \R_+^2.
\end{equation}

\smallskip  For the case (H3'), since
the flux function is continuous outside a negligible set $\mathcal
N$, then, for $x\in{\text{\menge R}}\backslash\mathcal N$,
$$
F(x,\tilde{\rho}(s,y,y))\stackrel{x\rightarrow
y}{\longrightarrow}F(y,\tilde{\rho}(s,y,y)).
$$
On the other hand, we have $F(y,\tilde{\rho}(s,y,y))
=F(x,\tilde{\rho}(s,y,x))$. Therefore, we have
$$
F(x,\tilde{\rho}(s,y,x))-F(x,\tilde{\rho}(s,y,y))\stackrel{x\rightarrow
y}{\longrightarrow}0,
$$
and \eqref{40} is a consequence of the fact that $F(x,\cdot)$ is a
one to one function. The case (H3) is clear from the definition of
$\tilde{\rho}(s,y,x)$ in \eqref{tilde-rho}.

\medskip
{\it Step 3.} With Steps 1--2, to achieve inequality \eqref{E}, it
suffices by choosing $\alpha=\beta(s,y)$ in \eqref{mves-2} to show
the following inequality:
\begin{eqnarray}\label{43}
&&\lim_{\tau,\omega\rightarrow 0} \int\partial_t J(t,x)\bar
H_{\tau}(t-s)H_{\omega}(x-y)
\langle\pi_{t,x}; |k-\tilde{\rho}(s,y,x)|\rangle \,dtdxdsdy\nonumber\\
&& +\lim_{\tau,\omega\rightarrow 0}\int\partial_x J(t,x)\bar
H_{\tau}(t-s)H_{\omega}(x-y)\langle \pi_{t,x};
{\text{sgn}}\left(k-\tilde{\rho}(s,y,x)\right)\big(F(x,k)-\beta(s,y)\big)\rangle\,
  dtdxdsdy\nonumber\\
&&+\lim_{\tau,\omega\rightarrow 0}\int J(0,x)\bar
H_{\tau}(-s)H_{\omega}(x-y)
|\rho_0(x)-\tilde{\rho}(s,y,x)|\,dxdsdy\,\,\geq\;0
\end{eqnarray}
for any test function $J\in C_0^{\infty}(\R_+^2)$ and verify that
$\tilde{\rho}(s,y,x)$ can be replaced by $\rho(t,x)$ in the limit as
$\tau, \omega\to 0$. Here the two families of functions $\bar
H_{\tau}$, $H_{\omega}\in C_0^{\infty}({\text{\menge R}})$ are
defined as
$$
\bar H_{\tau}(z)=\frac1{\tau}\bar H(\frac{z}{\tau})\quad\mbox{and}
\quad H_{\omega}(z)=\frac1{\omega}H(\frac{z}{\omega}) \qquad
\mbox{for}\, \tau,\omega>0,
$$
for a positive, compactly supported function ${H}\in
C_0^{\infty}({\text{\menge R}})$ and a positive function $\bar{H}\in
C_0^{\infty}({\text{\menge R}})$ with compact support in $(-2,-1)$
such that $\int_{\text{\menge R}}H(z)dz=1$ and $\int_{(-2,-1)}\bar
H(z)dz=1$.

This can be easily seen by first choosing the test function in
\eqref{mves-2} as $J(t,x)\bar{H}_{\tau}(t-s)H_{\omega}(x-y)\ge 0$
for fixed $(s,y)$ and then integrating with respect to $(s,y)$. We
now estimate the three terms of \eqref{43} in Steps 4--6,
respectively.

\medskip
{\it Step 4.} We show that, as $\tau,\omega \to 0$, the first term
converges to
$$
\int\partial_tJ(t,x)\,\langle \pi_{t,x};
|k-\rho(t,x)|\rangle\,dtdx.
$$
Observe that
\begin{eqnarray}\label{45}
&&\Big|\int\partial_tJ(t,x)\bar H_{\tau}(t-s)H_{\omega}(x-y)
\,\langle\pi_{t,x}; \big|k-\tilde{\rho}(s,y,x)\big|\rangle\, dtdxdsdy\nonumber\\
&&\,\,\,-\int\partial_tJ(t,x)\bar H_{\tau}(t-s)H_{\omega}(x-y)
\langle \pi_{t,x}; \big|k-\tilde{\rho}(s,y,y)\big|\rangle\,dtdxdsdy\Big|\nonumber\\
&&\leq\int |\partial_t J(t,x)|\bar H_{\tau}(t-s)\big(\int
H_{\omega}(x-y)\big|\tilde{\rho}(s,y,x)-\tilde{\rho}(s,y,y)\big|\,dx\big)\,
dtds dy \to 0 \quad\mbox{as}\,\, \omega\to 0,\qquad\quad
\end{eqnarray}
by the dominated convergence theorem and the fact that $\int
H_{\omega}(x-y)\big|\tilde{\rho}(s,y,x)-\tilde{\rho}(s,y,y)\big|\,dx\to
0$ as $\omega\to 0$ for a.e. $(s,y) \in \R_+^2$ since
$\tilde{\rho}(s,y,x)\stackrel{x\rightarrow y}{\longrightarrow}
\tilde{\rho}(s,y,y)=\rho(s,y)$ by \eqref{40}. Then, to find the
limit of the first term of \eqref{43}, it suffices to compute the
limit of
\begin{equation}\label{47}
\int\partial_t J(t,x)\bar
H_{\tau}(t-s)H_{\omega}(x-y)\langle\pi_{t,x};
\big|k-\rho(s,y)\big|\rangle\, dtdxdsdy.
\end{equation}
Thus, it suffices to show that $\rho(s,y)$ can be replaced by
$\rho(t,x)$ in \eqref{47}, i.e., as $\tau,\omega\to 0$,
\begin{eqnarray}\label{46}
&&\int\partial_t J(t,x)\bar H_{\tau}(t-s)H_{\omega}(x-y)
\big|\rho(t,x)-\rho(s,y)\big|\, dtdxdsdy\nonumber\\
&&=\int\partial_tJ(t,x)\bar H(-r)H(-z)\big|\rho(t,x)-\rho(t+\tau
r,x+\omega z)\big|\, dtdxdrdz\to 0.
\end{eqnarray}
This is guaranteed by the fact that
$$
\lim_{\tau,\omega\rightarrow0}\int\big|\rho(t,x)-\rho(t+\tau
r,x+\omega z)\big|\, dtdx=0,
$$
and the dominated convergence theorem since all the functions
involved are bounded. This implies that, in \eqref{47}, we can
indeed replace $\rho(s,y)$ by $\rho(t,x)$ to arrive at the result.

\smallskip
{\it Step 5.} We show that the second term of \eqref{43} converges
to
$$
\int\partial_x J(t,x)\langle\pi_{t,x}; {\text{sgn}}\left(k
-\rho(t,x)\right)\big(F(x,k)-\beta(t,x)\big)\rangle\, dtdx
\qquad\text{as}\,\, \tau,\omega \to 0.
$$

\smallskip The hypothesis (H2) on $F(x,\rho)$ implies
\begin{eqnarray*}
&&\Big|\,{\text{sgn}}\left(k-
 \tilde{\rho}(s,y,x)\right)\big(F(x,k)-\beta(s,y)\big)
 -{\text{sgn}}\left(k- \tilde{\rho}(s,y,y)\right)
 \big(F(x,k)-F(x, \tilde{\rho}(s,y,y))\big)\Big|\\
&&=\Big|\,{\text{sgn}}\left(k- \tilde{\rho}(s,y,x)\right)
  \big(F(x,k)-F(x, \tilde{\rho}(s,y,x))\big)
  -{\text{sgn}}\left(k-\rho(s,y)\right)\big(F(x,k)-F(x,\rho(s,y))\big)\Big|\\
&&\leq C| \tilde{\rho}(s,y,x)-\rho(s,y)|.
\end{eqnarray*}
Integrating the last expression with respect to $x$ against the
function $H_{\omega}(x-y)$ yields its convergence to zero by Step 2
as $\omega\rightarrow0$. Therefore, the limit of the
second term of \eqref{43} is the same as the limit of
$$
\int\!\!\partial_xJ(t,x)\bar
H_{\tau}(t-s)H_{\omega}(x-y)\langle\pi_{t,x}; {\text{sgn}}
\left(k\!-\!\rho(s,y)\right)\!\big(F(x,k)-F(x,\rho(s,y))\big)\rangle\,
dtdxdsdy,
$$
and it suffices to prove that, as $\tau,\omega\rightarrow0$,
$$
\begin{array}{ll}
\int\partial_x J(t,x)\bar
H_{\tau}(t-s)H_{\omega}(x-y)\langle\pi_{t,x};
\big|\,{\text{sgn}}\left(k-\rho(s,y)\right)
   \big(F(x,k)-F(x,\rho(s,y))\big)\\
\qquad\qquad \qquad\qquad\qquad\qquad\qquad\qquad
-{\text{sgn}}\left(k-\rho(t,x)\right)
    \big(F(x,k)-F(x,\rho(t,x))\big)\big|\rangle\,dtdxdsdy
\to 0.
\end{array}
$$
Using the Lipschitz property and fact \eqref{46}, we achieve the
result for the second term of \eqref{43}.

\medskip
{\it Step 6.} We now show that the third term of \eqref{43}
converges to zero if $\tau,\omega\rightarrow0$.
Note that
\begin{eqnarray*} \int J(0,x)\bar
H_{\tau}(-s)H_{\omega}(x-y)\Big||\rho_0(x)-\tilde{\rho}(s,y,x)|
-|\rho_0(x)-\tilde{\rho}(s,y,y)|\Big|\, dxdsdy\\
\leq\int J(0,x)\bar
H_{\tau}(-s)H_{\omega}(x-y)|\tilde{\rho}(s,y,x)-\tilde{\rho}(s,y,y)|\,
dxdsdy.
\end{eqnarray*}
For the same reason as in the first and the second term of
\eqref{43}, the right hand side converges to zero if $\tau,\omega\to
0$. We therefore next compute the limit as $\tau,\omega\rightarrow0$
of
\begin{equation*}\label{49}
\int J(0,x)\bar H_{\tau}(-s)H_{\omega}(x-y)|\rho_0(x)-\rho(s,y)|\,
dxdsdy.
\end{equation*}
As before, $\lim_{\tau,\omega\rightarrow0} \int
J(0,x)H_{\tau}(-s)H_{\omega}(x-y)|\rho(s,x)-\rho(s,y)|\, dxdsdy=0$.
Therefore, the next goal is to compute the limit of
\begin{equation}\label{50}
\int J(0,x) \bar{H}_{\tau}(-s)H_{\omega}(x-y)|\rho_0(x)-\rho(s,x)|\,
dxdsdy =\int J(0,x)\bar H(-r)|\rho_0(x)-\rho(\tau r,x)|\, dxdr.
\end{equation}
Since all the functions are bounded and $\mbox{supp}\, H\subset
(-2,-1)$, by the dominated convergence theorem and property
\eqref{initial-data} of the unique entropy solution $\rho(t,x)$,
this converges to $0$ as $\tau\to 0$ and thereby \eqref{50}
converges to $0$.

\smallskip With Steps 3--6 and by \eqref{43}, we complete the
proof.
\end{proof}

\smallskip
Then Theorem 2.1 yields the $L^1$-contraction between the
measure-valued entropy solution $\pi_{t,x}$ and the unique entropy
solution $\rho(t,x)$ of \eqref{1.1}--\eqref{1.1a}.

{\theorem [$L^1$-contraction]\label{l1} The function
$\int\langle\pi_{t,x}; |k-\rho(t,x)|\rangle\, dx$ is non-increasing
in $t>0$, which implies $\pi_{t,x}(k)=\delta_{\rho(t,x)}(k)$ when
$\pi_{0,x}(k)=\delta_{\rho_0(x)}(k)$ for a.e. $x\in\R$. }

\begin{proof} In expression (\ref{E}), we choose the test function as the product
test function $J_j(t)H(x)$, with $J_j(t)$ converging to the
indicator function $\text{\menge 1}_{[t_1,t_2]}(t)$ as
$j\rightarrow\infty$ for $t_2>t_1\ge0$. Then (\ref{E}) is equal to
\begin{eqnarray}
\nonumber&&\int\lefteqn{H(x) \langle \pi_{t_1,x}(k);
|k-\rho(t_1,x)|\rangle \,dx
-\int H(x)\langle \pi_{t_2,x}(k); |k-\rho(t_2,x)|\rangle\, dx}\\
&&\quad +\int_{t_1}^{t_2}\int H'(x){ \langle\pi_{t,x}(k);
\text{sgn}}\left(k
-\rho(t,x)\right)\big(F(x,k)-\beta(t,x)\big)\rangle\, dxdt\geq0.
\qquad\quad \label{2.18}
\end{eqnarray}
In \eqref{2.18}, we choose
$$
H(x)=e^{-\gamma\sqrt{1+|x|^2}}\chi(\frac{x}{N}), \qquad \gamma, N>0,
$$
for $\chi\in C_0^\infty(-2,2)$ with   $\chi(x)=1$ when $x\in [-1,1]$
and $\chi(x)\ge 0$. Letting $N\to\infty$ first and $\gamma\to 0$
then yields
that, for any $t_2>t_1\ge 0$,
$$
\int\langle \pi_{t_2,x}; |k-\rho(t_2,x)|\rangle\, dx - \int\langle
\pi_{t_1,x}; |k-\rho(t_1,x)|\rangle\, dx\le 0.
$$
In particular, when $t_2=t>0, t_1\to 0$, then
$\pi_{0,x}(k)=\delta_{\rho_0(x)}(k)$ implies
$$
\int\langle \pi_{t,x}; |k-\rho(t,x)|\rangle dx\le 0
$$
so that $\pi_{t,x}(k)=\delta_{\rho(t,x)}(k)$ for any $t>0$.
\end{proof}

\section{Existence of entropy solutions}

In this section, we establish the existence of  entropy solutions
(\ref{1.1})--\eqref{1.1a} in the sense of Definition 2.1, as
required for the reduction of measure-valued entropy solutions.
More precisely, for each fixed $\varepsilon>0$, $\rho^{\varepsilon}$
denotes the unique Kruzkov solution of (\ref{1.1})--\eqref{1.1a} in
the sense \eqref{3.2a}, where the flux function depends smoothly on
the space variable $x$; then it is shown that the sequence
$\rho^{\varepsilon}$ converges to an entropy solution of
(\ref{1.1})--\eqref{1.1a}.

\subsection[\small{Existence of entropy solutions when $F$ is smooth}]
{Existence of entropy solutions when $F$ is smooth}\label{3.1}

Define $F_{\varepsilon}(x,\rho)$ the standard mollification of
$F(x,\rho)$ in $x\in\R$:
\begin{equation}\label{smoothlambda}
F_{\varepsilon}(x,\rho):=(F(\cdot,\rho)\ast\theta^{\varepsilon})(x)
\to F(x,\rho) \quad a.e. \qquad\mbox{as}\,\, \varepsilon\to 0,
\end{equation}
with $\theta^{\varepsilon}(x):=\theta(\frac{x}{\varepsilon}),\,
\theta(x)\geq0, \, {\text{supp}}\, \theta(x)\subset[-1,1]$, and
$\int_{-1}^1\theta(x)dx =1$.  For fixed $\varepsilon>0$, consider
the
following Cauchy problem:
\begin{equation}\label{pdesmooth}
\Bigg\{\begin{array}{ll}
\partial_t\rho +\partial_xF_{\varepsilon}(x,\rho)=0,\\
\rho|_{t=0}=\rho_{0}(x)\ge 0.
\end{array}
\end{equation}

Kruzkov's result in \cite{kruzkov} indicates that there exists a
unique solution $\rho^{\eps}$ of \eqref{pdesmooth} satisfying the
Kruzkov entropy inequality:
\begin{equation}\label{3.2a}
\partial_t|\rho^\eps(t,x)-c | + \partial_x
\big(sgn(\rho^\eps(t,x)-c)
(F_\eps(x,\rho^\eps(t,x))-F_\eps(x,c))\big)
+(sgn(\rho^\eps(t,x)-c)\partial_xF_\eps(x,c)\leq0
\end{equation}
in the sense of distributions. We now show that the entropy solution
$\rho^\varepsilon$ also satisfies \eqref{es-2}.

{\proposition \label{KrAp} Let $\rho^{\varepsilon}(t,x)$ be a solution
of the Cauchy problem \eqref{pdesmooth} satisfying the Kruzkov
entropy inequality \eqref{3.2a}.
Then $\rho^\varepsilon(t,x)$ also satisfies the entropy inequality
\eqref{es-2} with
steady-state solutions $m_\alpha^\pm=m_\alpha^{\varepsilon,\pm}$. }

\begin{proof} In \eqref{3.2a}, we choose the
constant $c=m_{\alpha}^{\eps,\pm}(y)$ for any $\alpha\in
[M_0,\infty)$ (or $\alpha\in (-\infty, M_0]$), integrate against a
test function
$J_{\omega}(t,x,y):=J(t,\frac{x+y}2)H_{\omega}(x-y)$ with
$H_{\omega}$ as in the proof of Theorem \ref{dip}, integrate by
parts in the term involving $H_\omega'(x-y)$ with respect to $dy$,
and observe that $ (\partial_x+\partial_y)J_{\omega}(t,x,y)
=\partial_xJ(t,\frac{x+y}{2})H_{\omega}(x-y)$ to obtain from
\eqref{3.2a} that
\begin{eqnarray}\label{70}
\nonumber &&\int
 \left|\rho^\eps(t,x)-m_{\alpha}^{\eps,\pm}(y)\right|H_{\omega}(x-y)\partial_tJ\, dtdxdy\\
\nonumber&&+\int {\text{sgn}}
 \left(\rho^\eps(t,x)-m_{\alpha}^{\eps,\pm}(y)\right)
 \big(F_{\varepsilon}(x,\rho^{\eps}(t,x))-F_{\varepsilon}(x,m_{\alpha}^{\eps,\pm}(y))\big)
   H_{\omega}(x-y) \partial_yJ\, dtdxdy\\
\nonumber &&-\int\,{\text{sgn}}
 \left(\rho^\eps(t,x)-m_{\alpha}^{\eps,\pm}(y)\right)
 \partial_yF_{\varepsilon}(x,m_{\alpha}^{\eps,\pm}(y))H_{\omega}(x-y)\, J\, dtdxdy\\
 \nonumber &&-\int\, {\text{sgn}}
  \left(\rho^\eps(t,x)-m_{\alpha}^{\eps,\pm}(y)\right)
   \partial_xF_{\varepsilon}(x,m_{\alpha}^{\eps}(y))H_{\omega}(x-y) J\, dtdxdy\\
&&+\int\,
\left|\rho^\eps(0,x)-m_{\alpha}^{\eps,\pm}(y)\right|H_{\omega}(x-y)J(0,\frac{x+y}2)\,dxdy
\geq 0,
\end{eqnarray}
where we have used that
$$
\int\,
\big(F_{\varepsilon}(x,\rho^{\eps}(t,x))-F_{\varepsilon}(x,m_{\alpha}^{\eps,\pm}(y))
\big)H_{\omega}(x-y) J\,
d_y\big(\text{sgn}(\rho^\eps(t,x)-m_\alpha^{\eps,\pm}(y))\big)dxdt
=0.
$$

As in the proof of Theorem \ref{dip}, we can replace
$\rho^\eps(t,x)$ by $\rho^\eps(t,y)$ in the first term  as
$\omega\rightarrow0$ and replace $\rho^\eps(0,x)$ by
$\rho^\eps(0,y)$ in the last term as $\omega\rightarrow0$.

The second term is equal to
$$
\int {\text{sgn}}
 \left(\rho^\eps(t,x)-m_{\alpha}^{\eps,\pm}(y)\right)
 \big(F_{\varepsilon}(x,\rho^{\eps}(t,x))
 -F_{\varepsilon}(x,m_{\alpha}^{\eps,\pm}(y))\big)
  H_{\omega}(x-y) \partial_yJ\, dtdxdy.
$$
By the hypothesis (H3) on the flux function,
$F_{\varepsilon}(\cdot,\rho)$ is a Lipschitz function from
$(-\infty, \rho_m]$ and $[\rho_m,\infty)$ to $[M_0,\infty)$ (or
$(-\infty, M_0]$), which implies
$$
\left|{\text{sgn}}\left(\rho^\eps(t,x)-\rho^\eps(t,y)\right)
\big(F(x,\rho^\eps(t,x))-F(x,\rho^\eps(t,y))\big)\right| \leq
C|\rho^\eps(t,x)-\rho^\eps(t,y)|.
$$
One can show in a similar way as in the first term that
$$
\lim_{\omega\rightarrow0}
\int\left|\rho^\eps(t,x)-\rho^\eps(t,y)\right|H_{\omega}(x-y)
\partial_yJ\,dtdxdy=0.
$$
This means that, in the second term of \eqref{70}, one can replace
$F(x,\rho^\eps(t,x))$ by $F(x,\rho^\eps(t,y))$. Since $F_{\eps}$ is
also a smooth function with respect to the first variable, the
second term converges to
$$
\int\,{\text{sgn}}\left(\rho^\eps(t,y)-m_{\alpha}^{\eps,\pm}(y)\right)
\big(F_{\eps}(y,\rho^\eps(t,y))-\alpha\big)\partial_yJ(t,y)\, dtdy.
$$

In the third and fourth term in \eqref{70}, for $z\in{\text{\menge
R}}$, we have
\begin{eqnarray*}
&&\lim_{\omega\rightarrow0}\int
(\partial_x+\partial_y)F_{\eps}(x,m_{\alpha}^{\eps,\pm}(y)) \,
 \omega\,
 H_{\omega}(\omega z)J(t,y+\frac12\omega z)\, dz\\
&&=\lim_{\omega\rightarrow0}\int
   (\partial_x+\partial_y)F_{\eps}(x,m_{\alpha}^{\eps,\pm}(y))\,
    H(z)J(t,y+\frac12\omega z)\, dz\\
&&
  =J(t,y)\partial_y F_{\eps}(y,m_{\alpha}^{\eps,\pm}(y)\big)
=J(t,y)\partial_y\alpha =0.
\end{eqnarray*}

With these results, as $\omega{\rightarrow0}$, inequality \eqref{70}
becomes \eqref{es-2} for
$F_{\eps}(x,\rho)=(F(\cdot,\rho)\ast\theta^{\varepsilon})(x)$ with
steady-state solutions $m_\alpha^\pm=m_\alpha^{\varepsilon,\pm}$.
\end{proof}
Thus we conclude the existence of an entropy solution
$\rho_{\eps}(t,x)$ in the sense of Definition 2.1 for each
$F_{\eps}$ with fixed $\eps>0$.

\remark Notice that the sequence of approximate entropy solutions
converges to a measure-valued entropy solution as $\eps\to 0$:
First, since $\rho_0\in L^{\infty}$,
we find that, for $\alpha$ big enough,
$$
m_{\alpha}^{\eps,-}(x)\leq\rho_0(x)\leq m_{\alpha}^{\eps,+}(x)
\qquad\mbox{for all}\,\, x \in\R.
$$
{}From \cite{ap}, it then follows that
$$
m_{\alpha}^{\eps,-}(x)\leq\rho^{\eps}(t,x)\leq
m_{\alpha}^{\eps,+}(x),
$$
which implies the uniform boundedness  of $\rho^\varepsilon(t,x)$ in
$\varepsilon$ since $m_{\alpha}^{\eps,\pm}(x)$ are uniformly bounded
in $\varepsilon$.
Then there exists a compactly supported family of probability
measures $\pi_{t,x}$ on ${\text{\menge R}}$ (i.e. Young measures;
see Tartar \cite{tartar}) and a subsequence (still denoted by)
$\rho^\varepsilon(t,x)$ such that, for any continuous function
$f(\rho)$,
\begin{eqnarray}\label{yt}
f(\rho^{\varepsilon}(t,x)) \,\stackrel{*}{\rightharpoonup}\,
\langle\pi_{t,x}, f(k)\rangle \qquad \text{as} \,\,\varepsilon\to 0.
\end{eqnarray}
On the other hand, by Section \ref{3.1}, the sequence
$\rho^{\varepsilon}(t,x)$ satisfies the entropy inequality
\eqref{es-2} for $F_{\eps}(x,\rho)$ and the steady-state solutions
$m_\alpha^\pm=m_\alpha^{\eps,\pm}$. In particular, we use (\ref{yt})
and the definition of the sequence $F_{\varepsilon}(x,\rho)$ in
\eqref{smoothlambda} to conclude that, as $\varepsilon\to 0$, the
compactly supported family of probability measures $\pi_{t,x}$
satisfies that, for any test function $J: \R_+^2\mapsto{\text{\menge
R}_+}$,
\begin{eqnarray}
&&\int \big(\langle\pi_{t,x};
\big|k-m_{\alpha}^\pm(x)\big|\rangle\partial_tJ
+\left\langle\pi_{t,x}; {\text{sgn}}\left(k-m_{\alpha}^\pm\right)
\big(F(x,k)-\alpha\big)\right\rangle\partial_xJ\big)\, dxdt
\nonumber\\
&&\quad+\int \big|\rho_0(x)-m_{\alpha}^\pm(x)\big|J(0,x)\, dx\geq0.
 \label{59}
\end{eqnarray}
Thus, $\pi_{t,x}$ is a measure-valued entropy solution of
\eqref{1.1}--\eqref{1.1a} with compact support for a.e.
$(t,x)\in\R_+^2$ in the sense of Definition 2.2.

%
%%%%%%%%%%%%%%%%%%%%%%%%%%%%%%%%%%%%%%%%%%%%%%%%%%%%%%%%%%%%%%%%%%%%
\subsection
%[Existence of entropy solutions when $F$ is discontinuous in $x$]
{Existence of entropy solutions when $F$ is discontinuous in $x$}
\label{3.2}
%%%%%%%%%%%%%%%%%%%%%%%%%%%%%%%%%%%%%%%%%%%%%%%%%%%%%%%%%%%%%%%%%%%%%

We are now ready to state the main theorem of this section.

{\theorem \label{existence} Let $F(x,\rho)$ be strictly convex or
concave in $\rho$ for a.e. $x\in\R$ and satisfy {\rm (H1)--(H3)}, or
let $F(x,\rho)$ satisfy {\rm (H1)--(H2)} and {\rm (H3')}. Let
$\rho_0(x)\in L^{\infty}$. Then the sequence of entropy solutions
$\rho^{\eps}$ of the Cauchy problem \eqref{pdesmooth} (in the sense
of Definition {\rm 2.1})  converges to the unique entropy solution
of the Cauchy problem \eqref{1.1}--\eqref{1.1a} in the sense of
Definition {\rm 2.1}.}

\begin{proof} We consider the two cases separately.

For the case {\rm (H1)--(H2)} and {\rm (H3')}, that is, the flux
function $F$ is monotone in $\rho$, we apply the compactness
framework established in Section 2 to establish the convergence. For
this case, the existence of entropy solutions has been established
in \cite {BaitiJenssen}. In Remark 3.1, we have shown that the limit
of the entropy solutions $\rho^\eps$ is determined by a
measure-valued entropy solution $\pi_{t,x}$. Then, by Theorems
2.1--2.2, $\pi_{t,x}$ is the Dirac measure concentrated on the
unique entropy solution $\rho(t,x)$ of (\ref{1.1})--\eqref{1.1a} in
the sense of Definition 2.1, which implies the whole sequence
converges.

For the case {\rm (H1)--(H3)}, since we have not established the
existence of an entropy solution, we employ the compensated
compactness method to establish the convergence of the entropy
solutions of the Cauchy problem \eqref{pdesmooth}, which also yields
the existence of a unique entropy solution of the Cauchy problem
\eqref{1.1}--\eqref{1.1a}.

From Remark 3.1, we have known that $\rho^{\eps}$ is uniformly
bounded in $L^{\infty}$ which implies that there exists a
subsequence $\rho^{\eps}$ converging weakly to a compactly supported
family of probability measures $\nu_{t,x}$ on ${\text{\menge R}}_+$
such that, for any function $f(\rho, t,x)$ that is continuous in
$\rho$ for a.e. $(t,x)$,
\begin{eqnarray}
f(\rho^{\eps}(t,x), t,x) \,\,\stackrel{*}{\rightharpoonup}\,\,
\langle\nu_{t,x}, f(k, t,x)\rangle \qquad \text{as}\,\,
\varepsilon\to 0.
\end{eqnarray}
In particular,
\begin{equation}\label{rho-star}
\rho^{\eps}(t,x) \,\stackrel{*}{\rightharpoonup}\, \langle\nu_{t,x},
k\rangle=:\rho(t,x) \in L^\infty.
\end{equation}
Our goal is to prove the strong convergence of $\rho^{\eps}(t,x)$ to
$\rho(t,x)$ a.e., equivalently, $\nu_{t,x}=\delta_{\rho(t,x)}$,
which implies that $\rho(t,x)$ is an entropy solution of
(\ref{1.1})--\eqref{1.1a}, that is, $\rho(t,x)$ satisfies the
entropy inequality in Definition \ref{weakap}.

By Section 3.1, we have known that the sequence $\rho^{\eps}$ exists
and satisfies
\begin{eqnarray*}
E^\varepsilon:=
\partial_t\left|\rho_{\eps}(t,x)-\hat\rho^{\eps}(s,y,x)\right|
+\partial_x\big({\text{sgn}}\left(\rho^{\eps}(t,x)-\hat\rho^{\eps}(s,y,x)\right)
\left(F_{\eps}(x,\rho^{\eps}(t,x))-\gamma(s,y)\right)\big)\leq 0
\end{eqnarray*}
in the sense of distributions, where
\begin{eqnarray*}
\hat\rho^{\eps}(s,y,x)
:=m^{+,\eps}_{\gamma(s,y)}(x){\text{sgn}}_+(\rho(s,y)-\rho_m(y))
+m^{-,\eps}_{\gamma(s,y)}(x){\text{sgn}}_-(\rho(s,y)-\rho_m(y)).
\end{eqnarray*}
Notice that $\gamma(s,y):=F(y, \rho(s,y))$ is independent of $\eps$.
Thus, for fixed $(s,y)$, we have the strong convergence of
$m^{\pm,\eps}_{\gamma(s,y)}(x)$ to a steady-state solution
$m_{\gamma(s,y)}^\pm (x)$ of \eqref{1.1}--\eqref{1.1a} as $\eps\to
0$. In particular,
$$
\|\hat\rho^{\eps}\|_{L^{\infty}}\leq M, \qquad \text{$M$ independent
of $\eps$};
$$
and, for a.e. $(s,y,x)\in\R^2_+\times\R$,
\begin{eqnarray*}
\hat\rho^{\eps}(s,y,x)\, \to\,\hat\rho(s,y,x) :=
m^{+}_{\gamma(s,y)}(x){\text{sgn}}_+(\rho(s,y)-\rho_m(y))
+m^{-}_{\gamma(s,y)}(x){\text{sgn}}_-(\rho(s,y)-\rho_m(y)),
\end{eqnarray*}
as $\eps\to 0$. By Schwartz's lemma, $E^\eps$ is a sequence of
measures; by Murat's lemma \cite{murat81}, $E^\eps$ is uniformly
bounded measure sequence in the measure space. This implies that
\begin{equation}\label{cc1}
 {\text{$E^\eps$ $\qquad$ is compact in }} W^{-1,p}_{loc}(\R_+^2)
 \qquad \text{for any}\,\, p\in (1,2).
\end{equation}
On the other hand, since the vector field sequence
$$
(\big|\rho^{\eps}(t,x)-m_{\gamma(s,y)}^{\pm,\eps}(x)\big|,\,\,
{\text{sgn}}\big(\rho^{\eps}(t,x)-m_{\gamma(s,y)}^{\pm,\eps}(x)\big)
\left(F_{\eps}(x,\rho^{\eps}(t,x))-\gamma(s,y)\right))
$$
is uniformly bounded in $\eps$ for any fixed $(s,y)$, it follows
that
\begin{equation}\label{cc2}
 {\text{$E^\eps$ $\qquad$ is bounded in }} W^{-1,\infty}_{loc}(\R_+^2).
\end{equation}
With \eqref{cc1}--\eqref{cc2}, we obtain by a compact interpolation
theorem in \cite{chen-a,dcl} that
\begin{equation}\label{h11}
E^\eps \qquad \text{is compact in}\,\, H^{-1}_{loc}(\R_+^2).
\end{equation}
On the other hand,
\begin{equation}\label{h12}
\partial_t\rho^{\eps}+\partial_xF_{\eps}(x,\rho^{\eps})=0\qquad \text{which is
automatically compact in}\,\, H^{-1}_{loc}(\R_+^2).
\end{equation}

Moreover, since $\hat{\rho}^\eps(s,y,x)$ strongly converges a.e.,
then we find that, as $\eps\to 0$,
\begin{equation}\label{entropy-weak-1}
\begin{array}{lcl}
\eta^{\eps}_1(\rho^{\eps},t,x,s,y)
&:=&\left|\rho^{\eps}(t,x)-\hat\rho^{\eps}(s,y,x)\right|\\
&\stackrel{*}{\rightharpoonup}&\langle\nu_{t,x}(k); |k-\hat\rho(s,y,x)|\rangle\\
&=:&\langle\nu_{t,x}; \eta_1(k,t,x,s,y)\rangle,\\
q^{\eps}_1(\rho^{\eps},t,x,s,y)
&:=&{\text{sgn}}\left(\rho^{\eps}(t,x)-\hat\rho^{\eps}(s,y,x)\right)
\left(F_{\eps}(x,\rho^{\eps})-\gamma(s,y)\right)\\
&\stackrel{*}{\rightharpoonup}&
\langle\nu_{t,x}(k); {\text{sgn}}\left(k-\hat\rho(s,y,x)\right)
\left(F(x,k)-\gamma(s,y)\right)\rangle\\
&=:&\langle \nu_{t,x}; q_1(k,t,x,s,y)\rangle,\\
\eta^{\eps}_2(\rho^{\eps}(t,x))&:=&\rho^{\eps}(t,x)\\
&\stackrel{*}{\rightharpoonup}& \langle\nu_{t,x}(k); k\rangle=\rho(t,x)\\
&=:&\langle\nu_{t,x}; \eta_2(k)\rangle,\\
q^{\eps}_2(\rho^{\eps}(t,x),x)&:=&F_{\eps}(x,\rho^{\eps})\\
&\stackrel{*}{\rightharpoonup}&\langle\nu_{t,x}(k); F(x,k)\rangle\\
&:=&\langle \nu_{t,x}; q_2(k,x)\rangle,
\end{array}
\end{equation}
and
\begin{equation}\label{entropy-weak-2}
\begin{vmatrix}\eta_1(\rho^\eps(t,x), s,y,x)&q_1(\rho^\eps(t,x),s,y,x)\\
\eta_2(\rho^\eps(t,x))&q_2(\rho^\eps(t,x),x)\end{vmatrix}
\,\,\stackrel{*}{\rightharpoonup}\,\,
\left\langle\nu_{t,x}; \begin{vmatrix}\eta_1(k,s,y,x)&q_1(k,s,y,x)\\
\eta_2(k)&q_2(k,x)\end{vmatrix}\right\rangle,
\end{equation}
where
\begin{eqnarray*}
&&(\eta_1(k,t,x,s,y),\, q_1(k,t,x,s,y)) =(|k-\hat\rho(s,y,x)|, \,\,
   {\text{sgn}}\left(k-\hat\rho(s,y,x)\right))
\left(F(x,k)-\gamma(s,y)\right)),\\
&&(\eta_2(k), q_2(k,x))=(k, \,\, F(x,k)).
\end{eqnarray*}
Together \eqref{h11}--\eqref{h12} with
\eqref{entropy-weak-1}--\eqref{entropy-weak-2}, we apply the
Div-Curl lemma (see Tartar \cite{tartar} and Murat \cite{murat}) to
obtain
\begin{eqnarray*}
\left\langle\nu_{t,x}; \begin{vmatrix}\eta_1(k,s,y,x)&q_1(k,s,y,x)\\
\eta_2(k)&q_2(k,x)\end{vmatrix}\right\rangle=
\begin{vmatrix}\left\langle\nu_{t,x};\eta_1(k,s,y,x)\right\rangle
&\left\langle\nu_{t,x}; q_1(k,s,y,x)\right\rangle\\
\left\langle\nu_{t,x}; \eta_2(k)\right\rangle
&\left\langle\nu_{t,x}; q_2(k,x)\right\rangle
\end{vmatrix}
\end{eqnarray*}
for all $(s,y),(t,x)\in \R\backslash\mathcal M$ with $\mathcal M$ a
set of measure zero in $\R_+^2$. Thus, we have
\begin{eqnarray*}
\!\!\!&&\!\!\!\!\!\!\left\langle\nu_{t,x};
|k-\hat\rho(s,y,x)|F(x,k)-k\,\,
{\text{sgn}}\left(k-\hat\rho(s,y,x)\right)\left(F(x,k)-\gamma(s,y)\right)\right\rangle\\
\!\!\!&&\!\!\!\!\!\!=\left\langle\nu_{t,x};
|k-\hat\rho(s,y,x)|\right\rangle \left\langle\nu_{t,x};
F(x,k)\right\rangle-\left\langle\nu_{t,x},k\right\rangle
\left\langle\nu_{t,x}; {\text{sgn}}\left(k-\hat\rho(s,y,x))\right)
\left(F(x,k)-\gamma(s,y)\right)\right\rangle.
\end{eqnarray*}
Equivalently, we have
\begin{eqnarray*}
&&\left\langle\nu_{t,x}; |k-\hat\rho(s,y,x)|
\left(F(x,k)-\left\langle\nu_{t,x}; F(x,k)\right)\right\rangle\right\rangle\\
&&-\left\langle\nu_{t,x}; (k-\rho(t,x)){\text{sgn}}
\left(k-\hat\rho(s,y,x)\right)\left(F(x,k)-F(y,\rho(s,y))\right)\right\rangle=0.
\end{eqnarray*}
Since this is true for all $(s,y)$ and $(t,x)$ except on a set
$\mathcal M$ of measure zero, we then choose $(s,y)=(t,x)$ for
$(t,x)\in\R\backslash\mathcal M$ to obtain
\begin{eqnarray*}
\nonumber&&\left\langle\nu_{t,x};
|k-\rho(t,x)|\left(F(x,k)-\left\langle
\nu_{t,x}; F(x,k)\right)\right\rangle\right\rangle\\
&&-\left\langle\nu_{t,x};
(k-\rho(t,x)){\text{sgn}}\left(k-\rho(t,x)\right)
\left(F(x,k)-F(x,\rho(t,x))\right)\right\rangle=0,
\end{eqnarray*}
that is,
\begin{equation}\label{divcurl}
\left\langle\nu_{t,x}; |k-\rho(t,x)|\right\rangle
\left(F(x,\rho(t,x))-\left\langle\nu_{t,x};
F(x,k)\right\rangle\right)=0.
\end{equation}

There are two possibilities:

When $\left\langle\nu_{t,x}; |k-\rho(t,x)|\right\rangle=0$, then we
have $\nu_{t,x}(k)=\delta_{\rho(t,x)}(k)$.

When $\left\langle\nu_{t,x}; F(x,k)\right\rangle-F(x,\rho(t,x))=0$,
we note that
\begin{eqnarray*}
\left\langle\nu_{t,x}; F(x,k)\right\rangle-F(x,\rho(t,x))&=&
\left\langle\nu_{t,x}; F(x,k)-F(x,\rho(t,x))\right\rangle\\
&=&\langle\nu_{t,x};  F_\rho(x,\rho)(k-\rho)
 +\frac12 F_{\rho\rho}(x,\xi)(k-\rho)^2\rangle\\
&=&F_\rho(x,\rho)\left\langle\nu_{t,x}; k-\rho\right\rangle
 +\frac{1}{2}\langle\nu_{t,x}; F_{\rho\rho}(x,\xi)(k-\rho)^2\rangle\\
&=&\frac12\langle\nu_{t,x}; \int_0^1 \theta F_{\rho\rho}(x,
\theta\rho+(1-\theta)k)\, d\theta \, (k-\rho)^2\rangle.
\end{eqnarray*}
Since $F(x,\rho)$ is strictly convex or concave in $\rho$, we
conclude
\begin{equation}\label{delta}
\nu_{t,x}(k)=\delta_{\rho(t,x)}(k)\qquad\text{for $(t,x)$ a.e.}
\end{equation}

Therefore, we have
$$
\rho^\eps(t,x) \to \rho(t,x) \qquad a.e. \,\,\,\text{as}\,\,\eps \to
0.
$$
Since the limit is unique via the uniqueness result in \cite{ap},
the whole sequence $\rho^\eps(t,x)$ strongly converges to
$\rho(t,x)$ a.e. It is easy to check that $\rho(t,x)$ is the unique
entropy solution of the Cauchy problem \eqref{1.1}--\eqref{1.1a} in
the sense of Definition 2.1.
\end{proof}

\remark{The conditions on the flux function $F(x,\rho)$ in Theorem
3.1 for the non-monotone case can be relaxed as follows:
 $F(x,\rho)$ satisfies (H1)--(H3) and  is convex or concave with
$$
\mathcal{L}^1\{\rho\, :\, F_{\rho\rho}(x,\rho)=0\}=0 \qquad\text{for
\, a.e.}\, \, x\in\R,
$$
where $\mathcal{L}^1$ is the one-dimensional Lebesgue measure. }

\section{Hydrodynamic Limit of a Zero Range Process with
Discontinuous Speed-Parameter}\label{section4}

In Section 2, we have established a compactness framework for
approximate solutions via the reduction of measure-valued entropy
solutions of \eqref{1.1}--\eqref{1.1a} in the sense of Definition
2.1.
In this section we focus on a microscopic particle system for a Zero
Range Process (ZRP) with discontinuous speed-parameter $\lambda(x)$.
We apply the compactness framework to show the hydrodynamic limit
for the particle system, when the distance between particles tend to
zero, to the unique entropy solution of the Cauchy problem
\begin{equation}\label{4.1a}
\partial_t\rho+\partial_x\left(\lambda(x)h(\rho)\right)=0
\end{equation}
and with initial data:
\begin{equation}\label{4.1b}
\rho|_{t=0}=\rho_0(x)\ge 0,
\end{equation}
where
$h(\rho)$ is a  monotone function of $\rho$, and $\lambda(x)$ is
continuous in $x\in\text{\menge R}$ with $0<\lambda_1\le
\lambda(x)\le \lambda_2<\infty$ for some constants $\lambda_1$ and
$\lambda_2$, except on a closed set $\mathcal N$ of measure zero.
Then $m_{\alpha}^+=m_{\alpha}^-:=m_{\alpha}$ for $\alpha\in
[0,\infty)$.

Rezakhanlou in \cite{rez1} first established the hydrodynamic limit
of the processus des misanthropes (PdM) with constant
speed-parameter. Covert-Rezakhanlou \cite{rez} provided a proof of
the hydrodynamic limit of a PdM with nonconstant continuous
speed-parameter $\lambda$. In both proofs, the most important step
is to show an entropy inequality at microscopic level, which then
implies the (macroscopic) Kruzkov entropy inequality, when the
distance between particles tends to zero, and thereby implies the
uniqueness of limit points. In this section, we generalize this to
the case when the speed-parameter $\lambda$ has jumps for the
attractive Zero Range Process (ZRP). In \S 4.1, we analyze some
properties of the ZRP. In \S 4.2, we prove the one-dimensional
microscopic entropy inequality letting $\eps=\eps(N)=N^{-\sigma},
\sigma\in (0,1)$, for a ZRP with discontinuous speed-parameter as
$N\to\infty$. In \S4.3, we show the existence of measure-valued
solutions via the microscopic entropy inequality and how inequality
\eqref{es-2} follows.

%%%%%%%%%%%%%%%%%%%%%%%%%%%%%%%%%%
\subsection{Some properties of the microscopic interacting particle system}
\label{4.1}
%%%%%%%%%%%%%%%%%%%%%%%%%%%%%%%%%%%%%%%%

We consider a system of particles with conserved total mass and
evolving on a one-dimensional lattice $\text{\menge Z}$ according to
a Markovian law. With the Euler scaling factor $N$, the microscopic
particle density is expected to converge to a deterministic limit as
$N\rightarrow\infty$, which is characterized by a solution of a
conservation law. Under the Euler scaling, $\frac1N$ represents the
distance between sites. Obviously we have two space scales: The
discrete lattice ${\text{\menge Z}}$ as embedded in ${\text{\menge
R}}$ with ``vertices" $\frac uN$ and $u\in{\text{\menge Z}}$. In
this way, the distances between particles tend to zero if $N$
increases to infinity. Sites of the microscopic scale  $\text{\menge
Z}$ are denoted by the letters $u,v$ and correspond to the points
$\frac{u}{N}$, $\frac{v}{N}$ in the macroscopic
 scale $\text{\menge R}$.
Points of the macroscopic space scale $\text{\menge R}$ are denoted
by the letters $x,y$ and correspond to the sites $[xN]$, $[yN]$ in
the microscopic space scale, where [z] is the integer part of z. We
denote by $\eta_t (u)$ the number of particles at time $t>0$ at site
$u$. Then the vector $\eta_t= ( \eta_t (u):u\in\text{\menge Z})$ is
called a configuration at time $t$  with configuration space
$\text{\menge N}^{\text{\menge Z}}$.

In general, the ZRP can be described as follows: Infinitely many
indistinguishable particles are distributed on a $1$-dimensional
lattice. Any site of the lattice may be occupied by a finite number
of particles. Associated to a given site $u$ there is an exponential
clock with rate $\lambda_{\eps}(\frac uN)g(\eta(u))$ depending on
the macroscopic spatial coordinates. Each time the clock rings on
the site $u$, one of the particles jumps to the site $v$ chosen with
probability $p(u,v)$. The elementary transition probabilities $p$:
$\text{\menge Z} \mapsto [0, 1]$ are supposed to be
\begin{enumerate}
\item translation invariant: $p(x, y)= p(0, y-x)=:p(y-x)$;
\item normalized: $\sum_{y} p(x, y)=1$, $p(x, x)=0$;
\item assumed to be of finite range: $p(x, y)=0$ for $\left|y-x\right|$
sufficiently large;
\item irreducible: $p(0, 1)> 0$.
\end{enumerate}
Without loss of generality, we assume that $\sum_z p(z)z=\gamma=1$;
otherwise, for $\gamma\ne 1$, we replace the function $h(\rho)$ by
$h(\rho)/\gamma$ in the following argument. The rate $g:\text{\menge
N}\rightarrow \text{\menge R}_+$ is a positive, nondecreasing
function with $g(0)=0$, $g(+\infty)=+\infty$, and
\begin{equation}\label{4.3g}
\frac{g(k)}{k^2}\to 0 \qquad\quad\text{as }\,\, k\to \infty.
\end{equation}

With this description, the Markov process $\eta_t$ is generated by
\begin{equation}\label{gen}
NL^N_{\eps}f(\eta)=N\sum_{u,v}\lambda_{\eps}(\frac
uN)g(\eta(u))p(v-u)(f(\eta^{u,v})-f({\eta})).
\end{equation}
Here $N$ comes from the Euler scaling factor speeding the generator,
thus $\eta_t$ denotes a configuration on which this speeded
generator $NL^N_{\eps}$ has acted for time $t$, and $\eta ^{u, v}$
represents the
 configuration $\eta $ where one particle jumped from $u$ to $v$:
$$
\eta ^{u,v}(w)= \Bigg{\{}\begin{array}{ll}
\eta (w) & \mbox{if}\,\, w\ne u,v, \\
\eta (u)-1 & \mbox{if}\,\, w=u, \\
\eta (v)+1 & \mbox{if}\,\, w=v.
\end{array}
$$
For any $\eps=\eps(N)>0$ and for any constant $\alpha\geq0$, we define a
product measure given by
\begin{equation}
\tilde\nu_\alpha^N(\eta):=\prod_u\frac
1{Z\big(\alpha/\lambda_{\eps}(\frac
uN)\big)}\frac{\alpha^{\eta(u)}}{(\lambda_{\eps}(\frac
uN))^{\eta(u)}g(\eta(u))!}:=\prod_u\tilde{\nu}_\alpha^N(\eta(u)),
\end{equation}
where
$Z$ is a partition function equal to
\begin{equation}\label{Z}
Z\big(\frac{\alpha}{\lambda_{\eps}(\frac
uN)}\big)=\sum_{n=0}^{\infty}\frac{\alpha^n}{\left(\lambda_{\eps}(\frac
uN)\right)^ng(n)!}.
\end{equation}
Then the expected value of the
occupation variable $\eta(u)$ is equal to
$$
E_{\tilde\nu_\alpha^N}[\eta(u)] =\frac{\alpha}{\lambda_{\eps}(\frac
uN)} \frac{Z^{\prime}\big(\frac{\alpha}{\lambda_{\eps}(\frac
uN)}\big)} {Z\big(\frac{\alpha}{\lambda_{\eps}(\frac uN)}\big)}
:=R\big(\frac{\alpha}{\lambda_{\eps}(\frac uN)}\big).
$$
Now let $h$ be the inverse function of $R$ to obtain
$$
h\big(R\big(\frac{\alpha}{\lambda_{\eps}(\frac uN)}\big)\big)
=\frac{\alpha}{\lambda_{\eps}(\frac uN)}\quad \Rightarrow\quad
\lambda_{\eps}(\frac
uN)h\left(E_{\tilde\nu_\alpha^N}[\eta(u)]\right)=\alpha \quad
\Leftrightarrow \quad E_{\tilde\nu_\alpha^N}[\eta(u)]=m_\alpha(\frac
uN),
$$
where $m_{\alpha}$ is a steady-state solution to
\begin{equation}
\partial_t\rho+\partial_x\left(\lambda_{\eps}(x)h(\rho)\right)=0.
\end{equation}
Furthermore, it follows that
$$
E_{\tilde\nu_\alpha^N}[g(\eta(u))]=h\left(m_{\alpha}(\frac
uN)\right).
$$
{}From now on, we set
\begin{equation}\label{invarient-measure}
\mu^N_{m_\alpha}(\eta)=\prod_u\nu_{m_{\alpha}(\frac
uN)}(\eta(u)):=\prod_u\tilde{\nu}_{\lambda_\varepsilon(\frac{u}{N})
h(m_\alpha(\frac{u}{n}))}^N(\eta(u)).
\end{equation}

The important attribute of the ZRP with nonconstant speed-parameter
is that the \textit{product} measure $\mu^N_{m_\alpha}(\eta)$ is invariant
under the generator $NL_{\eps}^N$, i.e.,
\begin{equation}
\int L^N_{\eps}(f(\eta))d\mu^N_{m_\alpha}(\eta)=0.
\end{equation}

As a reasonable initial distribution, we choose {\textit{the product
measure $\mu_0^N(\eta)$ associated to a bounded density profile}}
defined as follows: For a bounded density profile $\rho_0\ge 0$, the
probability that particles at time $t=0$ are distributed with
configuration $\eta$ is equal to
\begin{equation}\label{mu}
\mu_0^N(\eta):=\prod_u\frac
1{Z(h(\rho_{u,N})/\lambda_\varepsilon(\frac{u}{N}))}
\frac{(h(\rho_{u,N}))^{\eta(u)}}
{(\lambda_{\varepsilon}(\frac{u}{N}))^{\eta(u)} g(\eta(u))!},
\end{equation}
where $\rho_{u,N}\ge 0$ is a sequence satisfying $
\lim_{N\rightarrow\infty}\int|\rho_{[Nx],N}-\rho_0(x)|dx=0$ for
$[Nx]$ as the integer part of $Nx$. With this definition, we say
that a sequence of probability measures $\mu^N$ is associated to a
density profile $\rho\ge 0$ if
$$
\lim_{N\rightarrow\infty}\langle \mu^N(\eta)\, ;
\,\big|\frac1N\sum_u J(\frac uN)\eta(u)-\int
J(x)\rho(x)dx\big|\rangle=0\qquad\mbox{for every test
function}\,\,J.
$$
Furthermore, let
\begin{equation}\label{mut}
\mu_t^N=S_t^N\ast\mu_0^N,
\end{equation}
where $S_t^N=e^{tNL^N_{\eps}}$ is the semigroup corresponding to the
generator $NL^N_{\eps}$. Then the attractiveness for two initial
measures $\mu_{\rho_0}^{N}$ and $\mu_{\omega_0}^{N}$ with profiles
$\rho_t$ and $\omega_t$, respectively, implies that
$$
\mu_{\rho_0}^{N}\leq\mu_{\omega_0}^{N}\,\,\,\Rightarrow\,\,\,
\mu_{\rho_t}^{N}\leq\mu_{\omega_t}^{N}
$$
is satisfied by the assumption that $g$ is a nondecreasing function.
Moreover, it is easy to prove that $\mu_{\rho_0}\leq\mu_{\omega_0}$
if $\rho_0\leq\omega_0$. It then follows by attractiveness that, for
any constant $\alpha$ such that $m_{\alpha}(x)\geq\rho_0(x)$, we obtain that the
inequality $\mu_0^N\leq\mu^N_{m_\alpha}$ implies
\begin{eqnarray}\label{nu}
S_t^N\mu_0^N\leq S_t^N\mu^N_{m_\alpha} =\mu^N_{m_\alpha}.
\end{eqnarray}
Since our initial distribution has a bounded density profile, then
the density profile remains bounded at later time $t$.

\par The goal in proving the
hydrodynamic limit of a ZRP is that, if we start from a
configuration $\eta_0$ distributed with an initial measure $\mu_0^N$
associated to the bounded density profile $\rho_0$, then the
configuration $\eta_t$ at later time $t$ is distributed with the
measure $\mu_t^N$ defined by \eqref{mut} and having density profile
$\rho(t,\cdot)$, where $\rho$ is the solution of the Cauchy problem
\eqref{4.1a}--\eqref{4.1b} in the sense of Definition 2.1. In other
words, our main theorem in this section is the following.

{\theorem[Hydrodynamic limit of an attractive ZRP with discontinuous
speed-parameter]\label{mainthm} Let $\eta_t$ be an attractive ZRP
with \eqref{4.3g} initially distributed by the measure $\mu_0^N$
associated to a bounded density profile $\rho_0:
\R_+^2\mapsto{\text{\menge R}}_+$ as defined in \eqref{mu}. Let
$\eps=\eps(N)=N^{-\sigma}, \sigma\in (0,1)$. Then, at later time
$t$,
\begin{equation}\label{hyd}
\lim_{N\rightarrow\infty}\langle\mu_t^N(\eta); \,\Big|\frac1N\sum_u
J(\frac uN)\eta_t(u)-\int J(x)\rho(t,x)dx\Big|\rangle=0
\end{equation}
for any test function $J: \R_+^2\to \R$, where $\rho$ is the unique
solution of the Cauchy problem \eqref{4.1a}--\eqref{4.1b} in the
sense of Definition {\rm 2.1}.}

\smallskip To achieve this, we have to establish an entropy
inequality in the sense of Definition \ref{weakap} at microscopic
level. This will be done in \S 4.2 by using the scaling relation
$\eps=\eps(N)=N^{-\sigma}, \sigma\in (0,1)$. Associated to each
configuration $\eta_t$, we may define the empirical measure viewed
as a random measure on $\text{\menge R}$ by
\begin{equation}\label{empirical}
\chi^N_t(x):=\frac1N\sum_{u}\eta_t(u)\delta_{\frac uN}(x).
\end{equation}
Then
$\langle\chi^N_t(\cdot), J(\cdot)\rangle=\frac1N\sum_u J(\frac
uN)\eta_t(u)$,
and we can rewrite \eqref{hyd} by
\begin{equation}\label{hyd1}
\lim_{N\rightarrow\infty}\langle\mu_t^N(\eta); \, \big|\langle
\chi^N_t(\cdot), J(\cdot)\rangle-\int J(x)\rho(t,x)dx\big|\rangle=0.
\end{equation}

%%%%%%%%%%%%%%%%%%%%%%%%%%%%%%%%%%%%%%%%
%\medskip
\subsection{The entropy inequality at microscopic level\label{4.2}}
The following proposition is essential towards the hydrodynamic
limit.

{\proposition[Entropy inequality at microscopic level for
$\varepsilon=N^{-\sigma}$ with $\sigma\in (0,1)$ as
$N\to\infty$]\label{ap} Let $m_{\alpha}^{\eps}$ be the steady-state
solutions of \eqref{pdesmooth} as defined in \eqref{1.2} with
$F_\varepsilon(x,\rho)=\lambda_\varepsilon(x)h(\rho)$.
Let $\eta_t$ be the ZRP generated by $NL^N_\varepsilon$ defined by
\eqref{gen} and initially distributed by the measure $\mu_0^N$
defined by \eqref{mu}.  Let $\eta^l(u)$ be the average density of
particles in large microscopic boxes of size $2l+1$ and centered at
$u$:
$$
\eta^l(u):=\frac1{2l+1}\sum_{|u-v|\leq l}\eta(v).
$$
Then, for every test function $J:\R_+^2\to\R_+$,
\begin{eqnarray}\label{micrap}
&&\lim_{l\rightarrow\infty}\lim_{N\rightarrow\infty}
\mu_t^N\Big\{\int_0^t\frac1{N}\sum_u\Big(\partial_sJ(s,\frac uN)
\big|\eta_s^l(u)-m^{\eps}_{\alpha}(\frac uN)\big|
+\partial_xJ(s,\frac uN)\big|\lambda_\varepsilon(\frac
uN)h(\eta_s^l(u))-\alpha\big|\Big)ds\nonumber\\
&&\qquad\qquad\qquad\quad + \frac1{N}\sum_u J(0,\frac uN)
\big|\eta_0^l(u)-m^{\eps}_{\alpha}(\frac uN)\big|
\geq-\delta\Big\}=1.
\end{eqnarray}
}

\smallskip
Inequality \eqref{micrap} is the entropy inequality \eqref{es-2}
with $\rho$ replaced by the average density of particles in the
microscopic boxes of length $2l+1$.
\smallskip To prove the microscopic entropy inequality, we
consider the coupled process $(\eta_t,\xi_t)$ generated by
$N\bar{L}_{\eps}^N$, where $\bar{L}_{\eps}^N$ is defined by
\begin{eqnarray}\label{cgen}
\bar{L}_{\eps}^Nf(\eta,\xi) &&=\sum_{u,v}p(v-u)\lambda_{\eps}(\frac
uN)
 \min\{g(\eta(u)),g(\xi(u))\}\left(f(\eta^{u,v}\!\!,\xi^{u,v})\!-f(\eta,\xi)\right)
 \nonumber\\
\nonumber&&\,\,\,\,\,+\sum_{u,v}p(v-u)\lambda_{\eps}(\frac uN)
  \{g(\eta(u))-g(\xi(u))\}_+\left(f(\eta^{u,v},\xi)-f(\eta,\xi)\right)\\
&&\,\,\,\,\,+\sum_{u,v}p(v-u)\lambda_{\eps}(\frac
uN)\{g(\xi(u))-g(\eta(u))\}_+\left(f(\eta,\xi^{u,v})-f(\eta,\xi)\right).
\end{eqnarray}
Furthermore, denote the initial distribution of $(\eta_t,\xi_t)$ by
$\bar{\mu}_0^N=\mu_0^N\times\mu^{N}_{m^{\eps}_{\alpha}}$, where
$\mu_0^N$ is the initial measure with density profile $\rho_0$
defined by \eqref{mu} and $\mu^N_{m^{\eps}_{\alpha}}$ denotes the
invariant measure as defined in \eqref{invarient-measure}.

Then, to prove Proposition \ref{ap}, it suffices to prove the
following proposition.

{\proposition{\label{apsmooth1}} Let $(\eta_t,\xi_t)$ be the coupled
process, starting from $\bar{\mu}_0^N$, generated by
$N\bar{L}_{\eps}^N$ as defined by \eqref{cgen}. Let
$\bar\mu_t^N=\bar S_t^N\ast\bar{\mu}_0^N$, where $\bar S_t^N$ is the
semigroup corresponding to the generator $N\bar{L}_{\eps}^N$. Then,
for every test function $J:\R_+^2\to\R_+$ and every
$\varepsilon=N^{-\sigma}$ with $\sigma\in (0,1)$,
$$
\begin{aligned}
&\lim_{l\rightarrow\infty}\lim_{N\rightarrow\infty}
\bar\mu_t^N\,\bigg\{\!\int_0^T\!\!\frac1{N}\sum_u
\Big\{\partial_sJ(s,\frac uN)\left|\eta_s^l(u)-\xi_s^l(u)\right|
+\partial_xJ(s,\frac uN)\lambda_\varepsilon(\frac{u}{N})
\left|h(\eta_s^l(u))-h(\xi_s^l(u))\right|\!\Big\}ds\\
&\qquad\qquad\qquad\,\, + \frac1{N}\sum_u J(0,\frac uN)
\left|\eta_0^l(u)-\xi_0^l(u)\right|\geq\!-\delta\,\bigg\}\!=\!1.
\end{aligned}
$$
}

\medskip

Recall that a microscopic entropy inequality leading to the Kruzkov
entropy inequality has been proved in \cite{rez} for the process of
PdM with nonconstant but continuous speed-parameter
$\lambda_{\eps}$. Since  there does not exist an invariant product
measure for a PdM in general such that
$E_{\mu^N_{m^{\eps}_{\alpha}}}[\xi(u)]=m^{\eps}_{\alpha}(\frac uN)$,
to replace the process $\xi$  by the process
$m^{\eps}_{\alpha}(\frac uN)$, one has to apply the relative entropy
method of Yau \cite{yau}.

In our case of a space-dependent ZRP, the invariant product measure
is available so that we can approximate  the steady-state solution
$m^{\eps}_{\alpha}$  by a process $\xi$ distributed by the invariant
measure $\mu^N_{m^{\eps}_{\alpha}}$ for any $\alpha\in (0,\infty)$.
Then, Proposition \ref{ap} indeed directly follows from Proposition
4.2.

\subsection{Proof of Proposition \ref{apsmooth1}}
We split the proof in three steps.

{\it Step 1:} {\textit{Lower bound for the martingale.}} For a test
function $J$ with compact support in $\R_+^2$,
define by
$M_t^J$ the martingale vanishing at time $t=0$:
\begin{eqnarray*}
M_t^J&=&\frac1{N}\sum_uJ(t,\frac uN)\left|\eta_t(u)-\xi_t(u)\right|
 -\frac1{N}\sum_uJ(0,\frac uN)\left|\eta_0(u)-\xi_0(u)\right|\\
&&-\int_0^t(\partial_s+N\bar{L}_{\eps}^N)\big(\frac1{N}\sum_uJ(s,\frac
uN)\left|\eta_s(u)-\xi_s(u)\right|\big)ds.
\end{eqnarray*}
Since $J$ has compact support, then, for $t$ large enough,
\begin{eqnarray*}
M_t^J=-\frac1{N}\sum_uJ(0,\frac uN)\left|\eta_0(u)-\xi_0(u)\right|
-\int_0^t(\partial_s+N\bar{L}_{\eps}^N)\Big(\frac1{N}\sum_uJ(s,\frac
uN)\left|\eta_s(u)-\xi_s(u)\right|\Big)ds.
\end{eqnarray*}
We now calculate
\begin{eqnarray}\label{LN}
&&\bar{L}_{\eps}^N\left|\eta(u)-\xi(u)\right|\nonumber\\
&&=\sum_{v,w}p(w-v)\lambda_{\eps}(\frac vN)
\Big\{\min\{g(\eta(v)),g(\xi(v))\}
\big(\left|\eta^{v,w}(u)-\xi^{v,w}(u)\right|-\left|\eta(u)-\xi(u)\right|\big)\nonumber\\
&&\nonumber\qquad\qquad\qquad\qquad\quad\,\,\,\,\,+
  \{g(\eta(v))-g(\xi(v))\}_+\big(\left|\eta^{v,w}(u)-\xi(u)\right|
   -\left|\eta(u)-\xi(u)\right|\big)\\
&&\qquad\qquad\qquad\qquad\quad\,\,\,\,\,+\{g(\xi(v))-g(\eta(v))\}_+
  \big(\left|\eta(u)-\xi^{v,w}(u)\right|-\left|\eta(u)-\xi(u)\right|\big)\Big\}\nonumber\\
&&=\sum_v\big(1-G_{u,v}(\eta,\xi)\big)\Big(-p(v-u)\lambda_{\eps}(\frac
uN)
     \left|g(\eta(u))-g(\xi(u))\right|+p(u-v)\lambda_{\eps}(\frac vN)
     \left|g(\eta(v))-g(\xi(v))\right|\Big)\nonumber\\
&&\,\,\,-\sum_vG_{u,v}(\eta,\xi)\Big(p(v-u)\lambda_{\eps}(\frac uN)
       \left|g(\eta(u))-g(\xi(u))\right|+p(u-v)\lambda_{\eps}(\frac
vN)\left|g(\eta(v))-g(\xi(v))\right|\Big), \label{lN}
\end{eqnarray}
where $G_{u,v}$ is the indicator function that equals to $1$ if
$\eta$ and $\xi$ are not ordered, i.e.,
\begin{eqnarray*}
G_{u,v}(\eta,\xi)={\text{\menge
1}}\left\{\eta(u)<\xi(u);\eta(v)>\xi(v)\right\}+{\text{\menge
1}}\left\{\eta(u)>\xi(u);\eta(v)<\xi(v)\right\}.
\end{eqnarray*}

Notice that the second sum is nonpositive. Therefore, plugging in
the last expression in the martingale $M_t^J$ and then interchange
$u$ and $v$ in the last term, we can bound the martingale below by
\begin{eqnarray*}
\!\!\!\!&-&\!\!\!\!\!\frac1{N}\sum_uJ(0,\frac uN)\left|\eta_0(u)-\xi_0(u)\right|
-\int_0^t\frac1N\sum_u\partial_sJ(s,\frac uN)\left|\eta_s(u)-\xi_s(u)\right|ds\\
\!\!\!\!&+&\!\!\!\!\!\int_0^t\sum_{u,v}\big(J(s,\frac uN)-J(s,\frac
vN)\big)p(v-u)\big(1-G_{u,v}(\eta_s,\xi_s)\big)\,\lambda_{\eps}(\frac
uN)\left|g(\eta_s(u))-g(\xi_s(u))\right|\, ds.
\end{eqnarray*}
Since the transition probability $p$ is of finite range, i.e.
$p(z)=0$ if $|z|>r$ for some $r$, then
$$
\left(J(s,\frac uN)-J(s,\frac
vN)\right)p(v-u)=-\frac1N(v-u)p(v-u)\partial_xJ(s,\frac
uN)+O(\frac1{N^2}).
$$
With $v=u+y$, it then follows that the martingale is bounded below
by
$$
\begin{aligned}
&-\int_0^t\frac1N\sum\limits_u\Big\{\partial_sJ(s,\frac
uN)\big|\eta_s(u)-\xi_s(u)\big|\\
&\qquad\qquad\qquad\,\,  +\partial_xJ(s,\frac
uN)\lambda_{\eps}(\frac
  uN)\tau_u\big(\sum\limits_yyp(y)(1-G_{0,y})\big)
  \big|g(\eta_s(0))-g(\xi_s(0))\big|\Big\}ds\\
&\qquad \quad -\frac1{N}\sum\limits_uJ(0,\frac uN)
\left|\eta_0(u)-\xi_0(u)\right|
   +O(\frac1{N}).
\end{aligned}
$$

{\it Step 2:}  We  show
\begin{equation}\label{squaremart}
\lim_{N\rightarrow\infty}E_{\bar\mu_t^N}\big[\left(M_t^J\right)^2\big]=0.
\end{equation}

Recall that
$$
N_t^J:=(M_t^J)^2-\int_0^t\Big(N\bar
L^N_{\eps}(A^J(s,\eta,\xi))^2-2A^J(s,\eta,\xi)N\bar
L_{\eps}^N(A^J(s,\eta,\xi))\Big)\, ds
$$
is a martingale vanishing at time $t=0$, where $A^J$ is defined by
$$
A^J(t,\eta,\xi)=\frac1{N}\sum_uJ(t,\frac uN)|\eta_t(u)-\xi_t(u)|.
$$
Then, by definition, $E_{\bar\mu_s^N}\left[N_s^J\right]=0$ for all
$0\leq s\leq t$. Thus, it suffices to show that the expectation of
the integral term of $N_t^J$ converges to zero as $N\to \infty$. In
order to prove this, we first find that, by careful calculation,
\begin{eqnarray*}
&&N\bar L_{\eps}^N(A^J(s,\eta,\xi))^2-2NA^J(s,\eta,\xi)\bar L_{\eps}^N(A^J(s,\eta,\xi))\\
&&=\sum_{v,w}p(w-v)N\lambda_{\eps}(\frac vN)
 \Big\{\left|g(\eta_s(v))-g(\xi_s(v))\right|\frac1{N^{2}}\big(1-G_{v,w}(\eta_s,\xi_s)\big)
 \big(J(s,\frac wN)-J(s,\frac vN)\big)^2\\
&&\qquad\qquad\qquad\qquad\qquad\quad
  +\left|g(\xi_s(v))-g(\eta_s(v))\right|\frac1{N^{2}}G_{v,w}(\eta_s,\xi_s)
  \big(J(s,\frac vN)+J(s,\frac wN)\big)^2 \Big\}.
\end{eqnarray*}
Since $J$ is a smooth function, the first term of this expression is
less $\mathcal{O}(\frac{g(CN)}{N^2})$ for some constant $C$
depending on the total initial mass and therefore converges to zero
as $N\to\infty$ by \eqref{4.3g}. For the second term, we know that
$(J(s,\frac vN)+J(s,\frac wN))^2\leq4\left\|J\right\|_\infty^2$,
which implies
\begin{eqnarray*}
&&N\bar L_{\eps}^N(A^J (s,\eta,\xi))^2-2NA^J(s,\eta,\xi)\bar
L_{\eps}^N(A^J(s,\eta,\xi))\\
&&=\mathcal
O(\frac{g(CN)}{N^2})+\frac{4\left\|J\right\|_\infty^2}{N}\sum_{v,w}
G_{v,w}(\eta_s,\xi_s)p(w-v)\lambda_\eps(\frac
vN)\left|g(\xi_s(v))-g(\eta_s(v))\right|.
\end{eqnarray*}

Then, to conclude the proof of \eqref{squaremart}, it suffices to
show
\begin{equation}\label{53}
E_{\bar\mu_t^N}\Big[\int_0^t\big(\sum_{v,w}
G_{v,w}(\eta_s,\xi_s)p(w-v)\lambda_\eps(\frac
vN)\left|g(\xi_s(v))-g(\eta_s(v))\right|\big)\,ds \Big]=\mathcal
O(1).
\end{equation}
For this, we use the martingale $M_t^J$ vanishing at $0$ with
$J\equiv1$, that is,
\begin{eqnarray*}
M_t:=\frac1{N}\sum_u|\eta_t(u)-\xi_t(u)|-\frac1{N}\sum_u|\eta_0(u)-\xi_0(u)|
-\int_0^t\frac1{N}\sum_u N\bar L_{\eps}^N|\eta_s(u)-\xi_s(u)|ds.
\end{eqnarray*}
By \eqref{lN}, the integral term of the martingale is equal to
\begin{eqnarray*}
\int_0^t\frac2{N}\sum_{u,v}NG_{u,v}(\eta_s,\xi_s)p(v-u)\lambda_{\eps}(\frac
uN)\left|g\left(\eta_s(u)\right)-g\left(\xi_s(u)\right)\right|ds,
\end{eqnarray*}
by interchanging $u$ and $v$ in some terms. Then we find
\begin{eqnarray*}
&&E_{\bar\mu_t^N}\Big[\int_0^t2\sum_{u,v}G_{u,v}(\eta_s,\xi_s)p(v-u)
 \lambda_{\eps}(\frac uN)
 \left|g\left(\eta_s(u)\right)-g\left(\xi_s(u)\right)\right|ds\Big]\\
&&=\,E_{\bar\mu_t^N}\Big[\int_0^t\frac1{N}\sum_u|\eta_0(u)-\xi_0(u)|ds\Big]
 -E_{\bar\mu_t^N}\Big[\int_0^t\frac1{N}\sum_u|\eta_t(u)-\xi_t(u)|ds\Big]\\
&&\leq
E_{\bar\mu_t^N}\Big[\int_0^t\frac1{N}\sum_u|\eta_0(u)-\xi_0(u)|ds\Big].
\end{eqnarray*}
Since we assumed that both marginals of $\bar\mu_t^N$ are bounded,
\eqref{53} follows, which leads to \eqref{squaremart}.

With the result of Step 1 and \eqref{squaremart} and using the
Chebichev inequality, we obtain
\begin{eqnarray}\label{68}
\!\!\!\!\!\!\!\!\!\!\!\!\nonumber&\!\!\!\!\!\!&\!\!\!\!\!\!\!\!\!\!\!\!\bar\mu_t^N\!\!
\Big\{\frac1{N}\sum_uJ(0,\frac
uN)\left|\eta_0(u)-\xi_0(u)\right|
+\int_0^t\frac1N\sum_u\big\{\partial_sJ(s,\frac uN)\big|\eta_s(u)-\xi_s(u)\big|\\
\!\!\!\!\!\!\!\!\!\!\!\!\nonumber&\!\!\!\!&\!\!\!\!
+\partial_xJ(s,\frac uN)\lambda_{\eps}(\frac
uN)\tau_u\big(\sum_yyp(y)(1-G_{0,y})(\eta,\xi)\big)
\big|g(\eta_s(0))-g(\xi_s(0))\big|\big\}ds+O(\frac1{N})<\!-\delta\!\Big\}\\
\leq&\!\!\!\!\!\!\!\!\!\!\!\!\!\!\!\!&\;\bar\mu_t^N\left\{M_t^J\;>\;\delta\!\right\}
\;\leq\;\bar\mu_t^N\left\{\left|M_t^J\right|\;>\;\delta\right\}\;
\leq\;\frac{1}{\delta^2}E_{\bar\mu_t^N}\big[\left(M_t^J\right)^2\big],
\end{eqnarray}
which converges to $0$ as $N\rightarrow \infty$, for all $\delta>0$.

\medskip
{\it Step 3.} We next use the following summation by parts formula: {\it For any
bounded function $a$ of $\eta(\cdot)$ with $a(0)=0$ and for any
smooth test function $J:{\text{\menge R}}\rightarrow{\text{\menge
R}}$,
we obtain that, for any $L>0$,
\begin{equation}\label{sbp}
\frac1{N}\sum_{|u|\leq LN}J(\frac
uN)a(\eta(u))=\frac1{N}\frac1{(2l+1)}\sum_{|u|\leq LN}J(\frac
uN)\sum_{|u-v|\leq l} a(\eta(v))+\mathcal
O(\frac{l\,\|J\|_{Lip}}{N}).
\end{equation}}
Since we restrict $\varepsilon=N^{-\sigma}, \sigma\in(0,1)$, then
$\|\lambda_{\varepsilon}\|_{Lip}\le C/\varepsilon=CN^\sigma$ and
$\mathcal{O}(\frac{l\|\lambda_\eps\|_{Lip}}{N})
=\mathcal{O}(\frac{l}{N^{1-\sigma}})\to 0$ as $N\to\infty$ so that
we can use this summation by parts formula \eqref{sbp} to replace
inequality \eqref{68} by
\begin{eqnarray}\label{69}
&&\lim_{l\rightarrow\infty}\lim_{N\rightarrow\infty}
\bar\mu_t^N\Big\{\frac1N\sum_u J(0,\frac
uN)\frac1{2l+1}\sum_{|z-u|\leq l}
\left|\eta_0(z)-\xi_0(z)\right|\nonumber\\
&&\nonumber\qquad\qquad\qquad+\int_0^t\frac1N\sum_u\partial_sJ(s,\frac
uN)\frac1{2l+1}
\sum_{|z-u|\leq l}\left|\eta_s(z)-\xi_s(z)\right|ds\\
&&\nonumber\qquad\qquad\qquad
+\int_0^t\frac{1}{N}\sum_u\partial_xJ(s,\frac
uN)\lambda_{\eps}(\frac uN)\frac1{2l+1}\nonumber\\
&&\qquad\qquad\qquad\qquad\,\times \sum_{|z-u|\leq
l}\tau_z\big(\sum_yyp(y)(1-G_{0,y})(\eta_s,\xi_s)\big)\big|g(\eta_s(0))-g(\xi_s(0))\big|\,ds
<-\delta\Big\}=0.\qquad\quad
\end{eqnarray}
Notice that, in \eqref{69}, since $J$ is a positive function, by the
triangle inequality, we can remove the sum inside the absolute value
in the first line. Following the same argument as in \cite{rez,rez1}
(also \cite{coz}), since we first set $\eps=\frac1{N^{\sigma}}$,
independent of $\lambda_\eps(x)$, we can obtain the following {\it
one block estimates}:
\begin{eqnarray}\label{4.27a}
\lim_{l\rightarrow\infty}\lim_{N\rightarrow\infty}E_{\bar\mu_t^N}
\Big\{\int_0^t\frac1N\sum_u\Big|\frac1{2l+1}\sum_{|u-z|\leq
l}|\eta_s(z)-\xi_s(z)| -|\eta_s^l(u)-\xi_s^l(u)|\Big|ds\Big\}=0,
\end{eqnarray}
and
\begin{eqnarray}
&&\lim_{l\rightarrow\infty}\lim_{N\rightarrow\infty}E_{\bar\mu_t^N}
\Big\{\int_0^t\frac1N\sum_u\tau_u\Big|\frac1{2l+1}\sum_{|z|\leq
l}\tau_z\big(\sum_yyp(y)(1-G_{0,y})(\eta_s,\xi_s)\big)
\big|g(\eta_s(0))-g(\xi_s(0))\big|\qquad\qquad\nonumber\\
&&\qquad\qquad\qquad\qquad\qquad\qquad\quad
-\left|h(\eta_s^l(0))-h(\xi_s^l(0))\right|\Big|ds\Big\}=0.
\label{4.27b}
\end{eqnarray}
Combining \eqref{69} with \eqref{4.27a}--\eqref{4.27b}, we complete
the proof of Proposition 4.2.

%%%%%%%%%%%%%%%%%%%%%%%%%%%%%%%%%%%%%%%
\subsection{Existence of measure-valued entropy solutions}\label{4.3}
%%%%%%%%%%%%%%%%%%%%%%%%%%%%%%%%%%%%%%%%%%

In this section, we prove that Theorem \ref{ap} implies the
existence of a measure-valued entropy solution associated to the
configuration $\eta_t$. We recall the empirical measure
$\chi^N_t(x)$ associated to a configuration $\eta_t$ in
\eqref{empirical}.
We define the Young measures associated to $\eta_t$ as follows:
\begin{equation}\label{young}
\pi^{N,l}_t(x,k):=\frac1N\sum_u\delta_{\frac
uN}(x)\delta_{\eta_t^l(u)}(k),
\end{equation}
which implies $\langle \pi^{N,l}_t; J\rangle=\frac1N\sum_uJ(\frac
uN,\eta_t^l(u))$ for any $J\in C_0({\text{\menge
R}}\times{\text{\menge R}}_+)$. If $E$ is the configuration space,
then these two measures are finite positive measures on $E$ and, for
any $J\in C_0({\text{\menge R}})$, they are related by the formula
\begin{equation}
\langle \pi^{N,l}_t; k J(x)\rangle \approx \langle \chi^N_t(\cdot);
J(\cdot)\rangle.
\end{equation}
Notice that, since there are jumps, the probability measure
$\mu_t^N$ defined by \eqref{mut} must be defined on the Skorohod
space $D[(0,\infty),E]$, which is the space of right continuous
functions with left limits taking values in $E$.  Then, using the
one to one correspondence between the configuration $\eta_t$ and the
empirical measure $\chi^{N}_t(\cdot)$, the law of $\chi^N$ with
respect to $\mu_t^N$ will give us a probability measure $Q^N$ on the
Skorohod space $D[(0,\infty),\mathcal M_+({\text{\menge R}})]$, for
the space $\mathcal M_+({\text{\menge R}})$ of finite positive
measures on ${\text{\menge R}}$ endowed with the weak topology.

In the same way, we can associate a probability measure $\tilde
Q^{N,l}$ on the space $D[(0,\infty),\mathcal M_+(\R_+^2)]$.

\noindent With these definitions, we can state the main theorem of
this section as follows.

\noindent {\theorem[Law of large numbers for the Young
measures]\label{llny} Let $(\mu^N)_{N\geq1}$ be a sequence of
probability measures, as defined by \eqref{mu}, associated to a
bounded density profile $\rho_0: {\text{\menge
R}}\rightarrow{\text{\menge R}}_+$. Then the sequence $\tilde
Q^{N,l}$ converges, as $N\rightarrow\infty$ first and
$l\rightarrow\infty$ second, to the probability measure $\tilde Q$
concentrated on the measure-valued entropy solution $\pi_{t,x}$ in
the sense of Definition {\rm 2.2}. }

\begin{proof}
In order to be allowed to take the limit points $Q$ and $\tilde Q$
of $Q^N$ and $\tilde Q^{N,l}$ respectively, we must know that the
sequences are tight (weakly relatively compact). If $Q^{N,l}$ is
weakly relatively compact, we can take $\tilde Q^l$ as a limit point
if $N\rightarrow\infty$ for each $l$. Denote by $\tilde Q$  a limit
point of $\tilde Q^{N,l}$ if $N\to\infty$ first and $l \to \infty$
second. Therefore, the proof consists in two main steps: The first
is to show that $\tilde Q^{N,l}$ is weakly relatively compact and
the second is to show the uniqueness of limit points. The key point
in the proof is that these can be achieved independent of the choice
of mollification $\lambda_\eps$ of the discontinuous speed-parameter
$\lambda$ with our choice of the notion of measure-valued entropy
solutions.

These can be achieved by following exactly the standard argument in
\cite{rez,rez1,kipnis} since it requires only the uniform
boundedness of $\lambda_\eps$ in the proof. That is, we can conclude
the following: Let $\mu_t^N$ be a measure defined by \eqref{mut}.
Then
\begin{enumerate}
\item[(i)] The sequence $Q^N$
defined above is tight in $D[(0,\infty),\mathcal M_+({\text{\menge
R}})]$ and all its limit points $Q$ are concentrated on weakly
continuous paths $\chi(t,\cdot)$;

\item[(ii)] Similarly, the sequence
$\tilde Q^{N,l}$ is tight in $D[(0,\infty),\mathcal
M_+({\text{\menge R}}\times{\text{\menge R}}_+)]$ and all its limit
points $\tilde Q$ are concentrated on weakly continuous paths
$\pi(t,\cdot,\cdot)$;

\item[(iii)] For every $t\geq0$, $\pi(t,x,k):=\pi_t(x,k)$ is
absolutely continuous with respect to the Lebesgue measure on
${\text{\menge R}}$, $\tilde Q$ a.s.. That is,  $\tilde Q$ a.s.
\begin{equation}\label{unique-1}
\pi_t(x,k)=\pi_{t,x}(k)\otimes dx;
\end{equation}

\item[(iv)] For every $t\in[0,T]$, $\pi_{t,x}(k)$ is compactly
supported, that is, there exists $k_0>0$ such that
$$
\pi_{t,x}([0,k_0]^c)=0\qquad\forall\,
(t,x)\in[0,T]\times{\text{\menge R}}.
$$

\item[(v)] $\pi_{t,x}$ is a measure-valued entropy
solution in the sense of Definition {\rm 2.2} for any $\alpha\in
[M_0,\infty)$, i.e.,
\begin{equation}\label{4.38a}
\partial_t\left\langle\pi_{t,x}; |k-m_{\alpha}(x)|\right\rangle+\partial_x\left\langle\pi_{t,x};
\left|h(k)\lambda(x)-\alpha\right|\right\rangle\leq0
\end{equation}
on $\R_+^2$ in the sense of distributions for any $\alpha\in
[M_0,\infty)$ or $\alpha\in (-\infty, M_0]$.
\end{enumerate}

The last result follows from the entropy inequality at microscopic
level in Theorem \ref{ap}. Indeed, in terms of the Young measures,
the expression \eqref{micrap} of Proposition 4.1:
\begin{eqnarray*}
&&\lim_{l\rightarrow\infty}\lim_{N\rightarrow\infty}\mu^N_t
\Big\{\int_0^{\infty}\frac1N\sum_{u}\big\{\partial_tH(t,\frac uN)
\big|\eta_t^l(u)-m_{\alpha}(\frac uN)\big|\\
&&\qquad\qquad\qquad\qquad\qquad\qquad+\partial_xH(t,\frac
uN)\big|\lambda(\frac
uN)h(\eta_t^l(u))-\alpha\big|\big\}dt\geq-\delta\Big\}=1
\end{eqnarray*}
can be restated as
\begin{eqnarray*}
&&\lim_{l\rightarrow\infty}\lim_{N\rightarrow\infty}\tilde Q^{N,l}
\Big\{\int_0^T
\big(\left\langle\pi_{t}(x,k);|k-m_{\alpha}(x)|\partial_tH(t,x)\right\rangle\\
&&\qquad\qquad\qquad\qquad\qquad+\left\langle\pi_{t}(x,k);\left|\lambda(x)h(k)-\alpha\right|
\partial_xH(t,x)\right\rangle\big) dt\geq-\delta\Big\}=1,
\end{eqnarray*}
for every smooth function $H:(0,T)\times {\text{\menge
R}}\rightarrow{\text{\menge R}}_+$ with compact support, any
$\alpha\in[M_0,\infty)$ or $\alpha\in (-\infty, M_0]$, and any
$\delta>0$. Since $\tilde Q$ is a weak limit point concentrated on
absolutely continuous measures and since we already proved that
$\pi_{t,x}$ is concentrated on a compact set (and therefore the
integrand is a bounded function), we obtain from the last expression
that
\begin{eqnarray*}
\tilde Q\Big\{\int_0^T\int\Big( \left\langle\pi_{t,x};
|k-m_{\alpha}(x)|\right\rangle\partial_tH(t,x)
+\left\langle\pi_{t,x};\left|\lambda(x)h(k)-\alpha\right|
\right\rangle\partial_xH(t,x)\Big)\, dx dt\geq-\delta\Big\}=1.
\end{eqnarray*}
Letting $\delta\rightarrow0$, we have that $\tilde Q$ a.s.
\eqref{4.38a} holds on $(0,T)\times\R$ in the sense of distributions
for every $\alpha\in[0,\infty)$. This proves the uniqueness of
$\tilde Q$ and thereby concludes the proof of Proposition
\ref{llny}.
\end{proof}

Then Theorem \ref{mainthm} follows immediately from this result
since the measure-valued entropy solution reduces to the Dirac mass
concentrated on the unique entropy solution $\rho(t,x)$ of
\eqref{4.1a}--\eqref{4.1b} as we noticed in \S3.2.

\medskip \noindent {\bf Acknowledgments}. Gui-Qiang Chen's
research was supported in part by the National Science Foundation
under Grants DMS-0505473, DMS-0426172, DMS-0244473, and an Alexandre
von Humboldt Foundation Fellowship. The research of Nadine Even and
Christian Klingenberg was support in part by a German DAAD grant
(PPP USA 2005/2006). The first author would like to thank Professor
Willi J\"{a}ger for stimulating discussions and warm hospitality
during his visit at the University of Heidelberg (Germany).

\bibliographystyle{amsalpha}

\end{document}